\documentclass[reqno,11pt]{amsart}
\usepackage{amsfonts}
\usepackage{graphicx}
\usepackage{amsfonts,amsmath, amssymb}
\usepackage{bbm}
\usepackage{hyperref}
\usepackage{marginnote}
\usepackage{color}

\numberwithin{equation}{section}

\newcommand{\dis}{\displaystyle}

\newtheorem{theorem}{Theorem}[section]
\newtheorem{corollary}[theorem]{Corollary}
\newtheorem{lemma}[theorem]{Lemma}

\newtheorem{remark}[theorem]{Remark}

\def\v{\varepsilon}

\def\t{\theta}

\def\n{\nu}
\def\m{\mu}
\def\a{\alpha}

\def\d{\delta}
\def\l{\lambda}

\def\r{\rho}
\def\s{\sigma}

\def\o{\omega}
\def\di{\displaystyle}
\def\i{\infty}
\def\f{\frac}
\def\pa{\partial}

\def\ltv{L^2_{v}}

\def\liv{L^\infty_{v}}
\def\w{\bar{w}}
\def\L{L^{-1}}
\def\fep{f^{\varepsilon}}
\def\hep{h^{\varepsilon}}
\def\Fep{F^{\varepsilon}}
\def\FRep{F^{\varepsilon}_{R}}
\def\ep{\varepsilon}

\def\pt{\partial_{t}}
\def\nax{\nabla_{x}}
\def\dis{\displaystyle}
\def\la{\langle}
\def\ra{\rangle}
\def\intr{\int_{\mathbb{R}^{3}}}
\def\ints{\int_{\mathbb{S}^{2}}}
\def\smum{\sqrt{\mu_{M}}}

\usepackage{tikz}

\begin{document}
	
	\title[Hydrodynamic limit of Boltzmann equation]{Hydrodynamic limit of the Boltzmann equation to the planar rarefaction wave in three dimensional space}

	\author[G.F. Wang]{Guanfa Wang}
	\address[G.F. Wang]{South China Normal University, Zhong Shan Avenue West 55, Tianhe District, Guangzhou 510631, China}
	\email{wgf@amss.ac.cn}

	\author[Y. Wang]{Yong Wang}
	\address[Y. Wang]{Institute of Applied Mathematics, Academy of Mathematics and Systems Science, Chinese Academy of Sciences, Beijing 100190, China, and University of Chinese Academy of Sciences}
	\email{yongwang@amss.ac.cn}
	
	\author[J.W. Zhou]{Jiawei Zhou}
	\address[J.W. Zhou]{School of Mathematical Sciences, University of Chinese Academy of Sciences, Beijing 100049, China; Institute of Applied Mathematics, AMSS, CAS, Beijing 100190, P.R.~China}
	\email{zhoujiawei15@mails.ucas.ac.cn}
	
\begin{abstract}
In this paper, we establish the global in time hydrodynamic limit of Boltzmann equation to the planar rarefaction wave of compressible Euler system in three dimensional space $x\in\mathbb{R}^3$ for general collision kernels.  Our approch is based on a generalized Hilbert expansion, and  a recent $L^2-L^\infty$ framework.  In particular, we improve the $L^2$-estimate to be a localized version because the planar rarefaction wave is indeed a one-dimensional wave which makes the source terms to be not integrable in the $L^2$ energy estimate of  three dimensional problem. We also point out that the wave strength of rarefaction may be large. 
\end{abstract}

\date{\today}
\maketitle
	
\setcounter{tocdepth}{1}
\tableofcontents
	
\thispagestyle{empty}

	
	\section{Introduction and main results}
	
\subsection{Introduction}	In this paper, we consider the Boltzmann equation
	\begin{equation}\label{1.1}
		F^\v_t+v\cdot\nabla_x F^\v=\frac1{\v}Q(F^\v,F^\v),
	\end{equation}
	where $F^\v(t,x,v)\geq 0$ is the density distribution function for the gas particles with position $x\in\mathbb{R}^3$  and velocity $v\in\mathbb{R}^3$ at  time $t>0$, and $\v>0$ is  Knudsen number which is proportional to the mean free path. The Boltzmann collision term $Q(F_1,F_2)$ on the right is defined in terms of the following bilinear form
	\begin{align}\label{1.2}
		Q(F_1,F_2)&\equiv\int_{\mathbb{R}^3}\int_{\mathbb{S}^2} B(v-u,\t)F_1(u')F_2(v')\,{d\omega du}
		-\int_{\mathbb{R}^3}\int_{\mathbb{S}^2} B(v-u,\t)F_1(u)F_2(v)\,{d\omega du}\nonumber\\
		&:=Q_+(F_1,F_2)-Q_-(F_1,F_2),
	\end{align}
	where the relationship between the post-collision velocity $(v',u')$ of two particles with the pre-collision velocity $(v,u)$ is given by
	\begin{equation*}
		u'=u+[(v-u)\cdot\omega]\omega,\quad v'=v-[(v-u)\cdot\omega]\omega,
	\end{equation*}
	for $\omega\in \mathbb{S}^2$, which can be determined by conservation laws of momentum and energy
	\begin{equation*}
		u'+v'=u+v,\quad |u'|^2+|v'|^2=|u|^2+|v|^2.
	\end{equation*}
The Boltzmann collision kernel $B=B(v-u,\theta)$ in \eqref{1.2} depends only on $|v-u|$ and $\theta$ with $\cos\theta=(v-u)\cdot \omega/|v-u|$.   Throughout this paper,  we  consider both the hard and soft potentials under the Grad's angular cut-off assumption, for instance,
	\begin{equation}\label{1.4}
		B(v-u,\t)=|v-u|^{\gamma}b(\t),
	\end{equation}
	with
	\begin{equation*}
		-3<\gamma\leq 1,\quad 0\leq b(\t)\lesssim|\cos\t|.
	\end{equation*}


It was shown since its derivation that the Boltzmann equation is closely related to the fluid dynamical systems for both compressible and incompressible flows, see the founding work of Maxwell \cite{Maxwell} and Boltzmann \cite{Boltzmann}.  For this, Hilbert proposed a systematic formal expansion in 1912, and Enskog and Chapman independently proposed another formal expansion in 1916 and 1917, respectively. Either the Hilbert expansion or the Chapman-Enskog expansion yield the compressible Euler equations in the leading order with respect to $\v$, and the compressible Navier-Stokes equations in the subsequent orders. It is a challenging problem to rigorously justify these formal approximation, that is, hydrodynamic limits. In fact, the purpose of Hilbert's sixth problem \cite{Hilbert} is to establish  the laws of motion of continua from the Boltzmann equation in mathematical standpoint.\vspace{1.5mm}

We review some previous works on the hydrodynamic limits of Boltzmann equation. For the case when the compressible Euler equations have smooth solutions, the hydrodynamic limits of the Boltzmann equation has been studied even in the case with an initial layer; cf. Caflisch \cite{Caflish}, Guo \cite{Guo Jang Jiang,Guo Jang Jiang-1},  Lachowicz \cite{Lachowicz}, Nishida \cite{Nishida}, and Ukai-Asano \cite{Ukai-Asano}.  However, as is well known, solutions of the compressible Euler equations in general develop singularities, such as shock waves. The Riemann problem was first formulated and studied by Riemann in the 1860s when he studied one-dimensional isentropic gas dynamics with initial data being two constant states. The Riemann solution turns out to be fundamental in the theory of hyperbolic conservation laws because it not only captures the local and global behavior of solutions but also fully represents the effect of nonlinearity in the structure of the solutions. It is now well known that for the  compressible Euler equations, there are three basic wave patterns, that is, shock wave, rarefaction wave, and contact discontinuity. These three types of waves have essential differences: shock is compressive, rarefaction is expansive, and contact discontinuity has some diffusive structure. Therefore, it is a natural problem to verify the hydrodynamic limit from the Boltzmann equation to the Euler equations with basic wave patterns.  For the one dimensional Boltzmann equation with slab symmetry, Yu \cite{Yu} proved the validity of hydrodynamic limit when the solution of the compressible Euler equations contains only noninteracting shocks;  Xin-Zeng \cite{Xin-Zeng} proved the case for rarefaction wave; the hydrodynamic limit to the contact discontinuity was proved by Huang-Wang-Yang \cite{Huang-Wang-Yang}. For superposition of different types of waves, we refer to the work \cite{Huang-Jiang-Wang,Huang-Wang-Yang-1,Huang-Wang-Wang-Yang}. In particular, Huang-Wang-Wang-Yang \cite{Huang-Wang-Wang-Yang} justify hydrodynamic limit in the setting of a Riemann solution that contains the generic superposition of shock, rarefaction wave, and contact discontinuity by introducing hyperbolic waves with different solution backgrounds to capture the extra masses carried by the hyperbolic approximation of the rarefaction wave and the diffusion approximation of contact discontinuity.\vspace{1.5mm}

For the case of  incompressible flows, the program was initiated by Bardos, Golse, and Levermore \cite{Bardos,Bardos-2} to justify the global weak solution of incompressible flows in the frame work of global renormalized solution of DiPerna-Lions \cite{Diperna-Lions}. In particular, Golse and Saint-Raymond \cite{Golse-Saint-Raymond} proved that the limits of the DiPerna-Lions renormalized solutions of the Boltzmann equation are the Leray solutions to the incompressible Navier-Stokes equations. There are also many important progresses on this topic such as \cite{Bardos-Ukai,E-Guo-M,E-Guo-K-M,Guo2006,Jiang-Masmoudi,Masmoudi-Raymond} and the references therein, we will not go into details about the incompressible limits since we will concentrate on the compressible Euler limit in this paper.\vspace{1.5mm}

We remark that all the works of hydrodynamic limit to the wave patterns of compressible Euler equations mentioned above are concerned in one-dimensional case, i.e. $x\in\mathbb{R}$. To the best of our knowledge, the hydrodynamic limit of Boltzmann equation to the wave patterns of compressible Euler equations in three dimensional space still remains open. The goal of this paper is to justify the limiting process of the Boltzmann equation to the planar rarefaction wave solution of compressible Euler equations in three dimensional case. The main difficulty is that the planar wave is indeed a one-dimensional wave in three dimensional space, and hence it and its derivatives are not integrable in $\mathbb{R}^3$. Therefore it is hard to use the  one-dimensional energy method to resolve it. To remedy the difficulty, we shall use a generalized Hilbert expansion, and a recent $L^2$-$L^\infty$ method \cite{Guo2010,Guo Jang Jiang}. In particular, we improve the $L^2$ estimation to be a localized version since the background planar rarefaction wave and its derivatives are not integrable in $\mathbb{R}^3$.\vspace{1.5mm}

\subsection{Hilbert expansion} We consider the  Hilbert expansion of Boltzmann solution \eqref{1.1}  with the form
	\begin{equation}\nonumber
		\Fep=\sum_{n=0}^{5}\ep^{n}F_{n}+\ep^{3}F^{\v}_{R}.
	\end{equation}
	where $F_{0},...,F_{5}$ are the first six terms of the Hilbert expansion, independent of $\ep$, which solve the equations:
	\begin{equation}\label{F_k equations}
		\begin{array}{l}
			\dis \v^{-1}: \quad\quad\quad\quad\quad\quad\ 0=Q(F_{0},F_{0}),\\[2mm]
			\dis \v^0:\quad \{\pt+v\cdot\nax\}F_{0}=Q(F_{0},F_{1})+Q(F_{1},F_{0}),\\[2mm]
			\dis \v: \quad\{\pt+v\cdot\nax\}F_{1}=Q(F_{0},F_{2})+Q(F_{1},F_{1})+Q(F_{2},F_{0}),\\[2mm]
			\dis \quad\quad\quad\quad\quad\quad\quad\ \vdots\\
			\dis \v^{5}:\quad  \{\pt+v\cdot\nax\}F_{5}=Q(F_{0},F_{6})+Q(F_{6},F_{0})+\sum_{\substack{i+j=6\\1\leq i,j\leq 5}}Q(F_{i},F_{j}).
		\end{array}
	\end{equation}
The reminder equation for $F_R^{\v}$ is given by
	\begin{align}\label{Remainder of FEP}
&\pt F^{\ep}_{R}+v\cdot \nabla_x F^{\ep}_R-\f{1}{\ep}\{Q(F_0, \FRep)+Q(\FRep, F_0)\}\nonumber\\
&=\ep^2Q(\FRep, \FRep)+\sum_{i=1}^5\v^{i-1}\{Q(F_i,\FRep)+Q(\FRep, F_i)\}\nonumber\\
&\quad +\ep^2\Big\{\sum_{\substack{i+j\geq 6\\1\leq i,j\leq 5}}\ep^{i+j-6}Q(F_{i},F_{j})-\{\partial_t+v\cdot \nabla_x\} F_5\Big\}.
\end{align}
It follows from  $\eqref{F_k equations}_1$ and the celebrated H-theorem that $F_0$ should be a local Maxwellian
\begin{equation}\label{1.6}
F_0(t,x,v)\equiv \f{\r_0(t,x)}{[2\pi \t_0(t,x)]^{3/2}}\exp{\left\{-\f{|v-u_0(t,x)|^2}{2\t_0(t,x)}\right\}},
\end{equation}
where $\r_0(t,x)$, $u_0(t,x)=(u_0^1,u_0^2,u_0^3)(t,x)$ and $\t_0(t,x)$ are defined as
\begin{align}\nonumber
&\intr F_0 dv =\r_0,\quad \intr vF_0dv =\r_0 u_0, \quad \intr|v|^2F_0 dv =\r_0|u_0|^2+3\r_0\t_0,
\end{align}
which represent the macroscopic density, velocity and temperature, respectively.  Projecting the equation  $\eqref{F_k equations}_2$ onto $1$, $v$, $\f{|v|^2}{2}$, which are five collision invariants for the Boltzmann collision operator $Q(\cdot,\cdot)$, one obtains that $(\r_0, u_0, \t_0)$ satisfies the compressible Euler system
	\begin{equation}\label{1.7}
\begin{cases}
\dis \pt\r_0+\mbox{div} (\r_0u_0)=0,\\[2mm]
\dis \pt(\r_0u_0)+\mbox{div}(\r_0u_0\otimes u_0)+\nabla p=0,\\[2mm]
\dis \pt[\r_0(\f{3\t_0}{2})+\f{|u_0|^2}{2}]+\mbox{div} [\r_0 u_0(\f{3\t_0}{2}+\f{|u_0|^2}{2})]+\mbox{div}(pu_0)=0,
\end{cases}
\end{equation}
where $p=\r_0\t_0$ is the pressure function.

\subsection{Planar rarefaction wave}
In this article we shall consider the hydrodynamic limit of Boltzmann equation to the planar rarefaction wave solution of compressible Euler equations. We impose \eqref{1.7} with the following
 Riemann initial data
\begin{equation}\label{Riemann-data}
(\rho_0,u_0,\theta_0)(0,x)=
\begin{cases}
(\r_-,u_{-},\t_-),~~~x_1<0,\\[1mm]
(\r_+,u_{+},\t_+),~~~x_1>0,
\end{cases}
\end{equation}
where $u_{\pm}=(u_{\pm}^1,0,0)$ and $\r_\pm>0,\t_\pm>0, u_{\pm}^1$ are given constants.
To construct a Riemann solution of \eqref{1.7} and \eqref{Riemann-data}, we introduce the  Riemann problem for the one dimensional inviscid Burgers equation:
\begin{equation}\label{bur}
\begin{cases}
w_t+ww_{x_1}=0,\\
w(x_1,0)=\begin{cases}
w_-,&x_1<0,\\
w_+,&x_1>0.
\end{cases}
\end{cases}
\end{equation}
If $w_-<w_+$,  the Riemann problem $(\ref {bur})$ admits a
rarefaction wave solution $w^r(x_1, t) = w^r(\frac{x_1}{t})$ given by
\begin{equation}\nonumber
w^r(\frac{x_1}{t})=
\begin{cases}
w_-,&\frac{x_1}{t}\leq w_-,\\
\frac{x_1}{t},&w_-\leq \frac{x_1}{t}\leq w_+,\\
w_+,&\frac{x_1}{t}\geq w_+.
\end{cases}
\end{equation}
In this paper we consider only the 1-rarefaction wave without loss of generality, since the 3-rarefaction wave can be treated similarly. Hence we assume that
$(\r_-,u_{-}^1,\t_-)$ and $(\r_+,u_{+}^1,\t_+)$ was connected by 1-rarefaction wave for the one-dimensional compressible Euler equations, then the Riemann
problem \eqref{1.7}, (\ref{Riemann-data}) admits a planar rarefaction wave
solution $(\r^{r_1}, u^{r_1},\t^{r_1})(t,x_1)$ defined by
\begin{align}\label{1.12}
\begin{cases}
\di  s^{r_1}=s(\r^{r_1},\t^{r_1})=s_+,\\
\di w^r(\frac{x_1}{t})= \l_1(\r^{r_1}(t,x_1),u^{1r_1}(t,x_1),s_+)=u^{1r_1}-\sqrt{\f{5}{3}}(\r^{r_1})^{\f13}{\rm exp}(\f{s_+}{2}),\\
\di u^{1r_1}(t,x_1)+\sqrt{15}(\r^{r_1}(t,x_1))^{\f13}{\rm
	exp}(\f{s_+}{2})=u_{+}^1+\sqrt{15}\r_+^\f13{\rm exp}(\f{s_+}{2}),\\
\di  u^{2r_1}=u^{3r_1}=0,
\end{cases}
\end{align}
where $\lambda_1(\rho,u^1,s)$ is the first eigenvalue of Euler equations, and $\dis s=\ln\theta-\frac23\ln\rho$ is the entropy.
It is noted that  a planar rarefaction wave is indeed a one dimensional wave in three dimensional space, and the wave is independent of the variables $x_2$ and $x_3$.

Notice that the planar rarefaction wave solution constructed in \eqref{1.12} is only Lipschitz continuous at the edge of the rarefaction wave, and has singularity at  $t=0, x=(0,x_2,x_3)$ for any $(x_2,x_3)\in\mathbb{R}^2$. To construct the linear part of Hilbert expansion $F_i, i=1,\cdots, 5$, we need more regularity on the planar rarefaction wave. Similar to \cite{Xin-1993,Huang-Wang-Wang-Yang}, we construct a smooth approximate 1-rarefaction wave.  Hence we consider Burgers equation
\begin{align}\label{Burgers equation}
\begin{cases}
\di w_{t}+ww_{x_1}=0,\\
\di w( 0,x_1
)=w_\sigma(x_1)=w(\f{x_1}{\sigma})=\f{w_++w_-}{2}+\f{w_+-w_-}{2}\tanh\f{x_1}{\sigma}.
\end{cases}
\end{align}
where $\sigma>0$ is a parameter  to be determined later. The solution $w^R_\sigma(t,x_1)$ of the Burgers equation is given by
\begin{align}\nonumber
w^R_\sigma(t,x_1)=w_\sigma(x_0(t,x_1)),\qquad
x_1=x_0(t,x_1)+w_\sigma(x_0(t,x_1))t.
\end{align}
where $x_0(t,x_1)=X(0;t,x_1)$, and $X(s;t,x_1)$ is the characteristic line of Burgers equation
\begin{align}\nonumber
\frac{dX(s)}{ds}=w(s,X(s)),\quad X(t;t,x_1)=x_1.
\end{align}
Then the approximate planar rarefaction wave $(\r^{R_1}, u^{R_1}, \t^{R_1}) (t,x_1)$  is given by
\begin{align}\label{construct 1 of r u t}
\begin{cases}
\di  s^{R_1}=s(\r^{R_1},\t^{R_1})=s_+,\\
\di w_\sigma^R(t,x_1)= \l_1(\r^{R_1}(t,x_1),u^{1R_1}(t,x_1),s_+)=u^{1R_1}-\sqrt{\f{5}{3}}(\r^{R_1})^{\f13}{\rm exp}(\f{s_+}{2}),\\
\di u^{1R_1}(t,x_1)+\sqrt{15}(\r^{R_1}(t,x_1))^{\f13}{\rm
	exp}(\f{s_+}{2})=u_{+}^1+\sqrt{15}\r_+^\f13{\rm exp}(\f{s_+}{2}),\\
\di  u^{2R_1}=u^{3R_1}=0,
\end{cases}
\end{align}
It is direct to know that the smooth approximate planar 1-rarefaction waves $(\r^{R_1}, u^{R_1},\theta^{R_1})(t,x_1)$ also satisfies the compressible Euler equations \eqref{1.7}. From Lemma \ref{lem2.1} below, we have that
\begin{align}\label{1.16}
&\sup_{x\in\mathbb{R}}|(\r^{R_1},u^{1R_1},\theta^{R_1})(t,x_1)-(\r^{r_1},u^{1r_1},\theta^{r_1})(\f{x_1}t)|\nonumber\\
&\leq\f{C}{t}[\s\ln(1+t)+\s|\ln\s|]\rightarrow 0,\quad\mbox{as} \ \sigma\rightarrow0,
\end{align}
for any given time $t>0$. That means the smooth approximate planar 1-rarefaction wave $(\r^{R_1}, u^{R_1},\theta^{R_1})(t,x_1)$ approximate the planar 1-rarefaction wave solution $(\r^{r_1}, u^{r_1},\t^{r_1})(t,x_1)$ very well after the initial time.

\subsection{Main results} 
From now on, we denote
\begin{equation}\label{2.10}
(\rho_0, u_0, \theta_0)(t,x):=(\r^{R_1}, u^{R_1},\theta^{R_1})(t,x_1).
\end{equation}
Then it is noted that $(\rho_0, u_0, \theta_0)(t,x)$ is a smooth solution to the compressible Euler equations \eqref{1.7}. We also define
\begin{align}\label{1.18}
F_0\equiv\mu_{\sigma}(t,x,v):=\f{\r_0(t,x)}{[2\pi \t_0(t,x)]^{3/2}}\exp{\left\{-\f{|v-u_0(t,x)|^2}{2\t_0(t,x)}\right\}},
\end{align}
and
\begin{align}\label{1.18-1}
\mu(t,x_1,v):=\f{\r^{r_1}(t,x_1)}{[2\pi \t^{r_1}(t,x_1)]^{3/2}}\exp{\left\{-\f{|v-u^{r_1}(t,x_1)|^2}{2\t^{r_1}(t,x_1)}\right\}},
\end{align}
where $(\rho_0, u_0, \theta_0)(t,x)$ and $(\r^{r_1}, u^{r_1},\t^{r_1})(t,x_1)$ are  the ones defined in \eqref{2.10} and \eqref{1.12}, respectively. We point out that the solution $(\rho_0, u_0, \theta_0)(t,x)$ depends on the parameter $\sigma$ throughout this paper even though we do not write it down explicitly.

For later use we define the linearized collision operator $\mathbf{L}$ by
\begin{equation}
\mathbf{L}g=-\frac{1}{\sqrt{\mu_\sigma}}\Big\{Q(\mu_{\sigma},\sqrt{\mu_\sigma} g)+Q(\sqrt{\mu_\sigma} g,\mu_\sigma)\Big\},\nonumber
\end{equation}
and the nonlinear operator
\begin{equation*}
\Gamma(g_1,g_2)=\frac{1}{\sqrt{\mu_\sigma}} Q(\sqrt{\mu_\sigma}g_1,\sqrt{\mu_\sigma}g_2).
\end{equation*}
The null space $\mathcal{N}$ of $\mathbf{L}$ is generated by
\begin{align}\nonumber
\begin{split}
\chi_0(v)&\equiv\f{1}{\sqrt{\r_0}}\sqrt{\mu_\sigma},\\
\chi_i(v)&\equiv\f{v^i-u_0^i}{\sqrt{\r_0\t_0}}\sqrt{\mu_\sigma},\quad i=1,2,3,\\
\chi_4(v)&\equiv\f{1}{\sqrt{6\r_0}}\left\{\f{|v-u_0|^2}{\t_0}-3\right\}\sqrt{\mu_\sigma}.
\end{split}
\end{align}
It is easy to check that $\displaystyle\int_{\mathbb{R}^3}\chi_i\cdot \chi_j dv=\delta_{ij}$ for $0\leq i,j\leq 4$.
We also define the collision frequency $\nu$:
\begin{equation}\label{2.3}
\nu(t,x,v)\equiv\nu(\mu_{\sigma}):=\intr\ints B(v-u,\t)\mu_\sigma(u)d\o du.
\end{equation}
It is direct to know that
\begin{equation}\nonumber
\frac1{C} (1+|v|)^{\gamma}\leq \nu(t,x,v) \leq  C(1+|v|)^{\gamma},
\end{equation}
where the constant $C>0$ depends only on $\rho_{\pm}, \theta_{\pm}, u_{\pm}$, but  is independent of $\sigma$.
Let $\mathbf{P}g$ be the $L_{v}^{2}$ projection with respect to $[\chi_0,...,\chi_4]$. It is well-known that there exists a positive number $c_{0}>0$ such that for any function $g$
\begin{equation}\nonumber
\la \mathbf{L}g,g\ra\geq c_{0}\|\{\mathbf{I-P}\}g\|_{\nu}^{2},
\end{equation}
where the weighted $L^2$-norm $\|\cdot\|_{\nu}$ is defined as
\begin{equation}\nonumber
\|g\|_{\nu}^2:=\int_{\mathbb{R}^3_x\times\mathbb{R}^3_v}  g^2(x,v) \nu(v)dxdv.
\end{equation}
We point out that the constant $c_0>0$ is independent of $\sigma$ even though the projection $\mathbf{P}$ depends on  $\sigma$.

For each $i\geq 1$, we define the macroscopic and microscopic part of $\displaystyle\frac{F_i}{\sqrt{\mu_{\sigma}}}$ as
\begin{align}\label{1.25}
\frac{F_i}{\sqrt{\mu_{\sigma}}}
&=\mathbf{P}\left(\frac{F_i}{\sqrt{\mu_{\sigma}}}\right)+\{\mathbf{I-P}\}\left(\frac{F_i}{\sqrt{\mu_{\sigma}}}\right)\nonumber\\
&\equiv \left\{\frac{\rho_i}{\sqrt{\rho_0}} \chi_0+\sum_{j=1}^3\sqrt{\frac{\rho_0}{\theta_0}} u_i^j\cdot \chi_j+\sqrt{\frac{\rho_0}{6}}\frac{\theta_i}{\theta_0} \chi_4 \right\}+\{\mathbf{I-P}\}\left(\frac{F_i}{\sqrt{\mu_{\sigma}}}\right).
\end{align}
\begin{theorem}[Estimates on the linear terms]\label{thm1.2}
Let  $\sigma\in(0,1]$,  $(\rho_0, u_0, \theta_0)(t,x)$ be the smooth approximate planar rarefaction wave of Euler equations constructed in \eqref{2.10}, and $F_0$ defined in \eqref{1.18}. For each $i\geq1$, we assume the initial data of macroscopic part
\begin{equation}\label{1.18-2}
(\r_{i},u_{i},\t_{i})(0,x_1):=(\r_{i0},u_{i0},\t_{i0})(x_1)\in H^s(\mathbb{R}),
\end{equation}
where $s>0$ is some positive constant, and  $\|(\r_{i0},u_{i0},\t_{i0})\|_{H^s}$ is independent of $\sigma>0$. Then the linear problem \eqref{F_k equations} is well-posed. Furthermore, there exists positive constants $C_0, C_{i},C_{i,n}\geq1, i=1,\cdots 5, n=1,\cdots$   such that
\begin{align}
|F_i(t,x_1,v)|&\leq C_{i}(\s+t)^{C_{0}i}\s^{-C_{0}i}(1+|v|)^{3i+(i-1)\bar\gamma}\mu_\sigma,\label{3.26-1}\\
|\pa_\tau^nF_i(t,x_1,v)|&\leq C_{i,n}(\s+t)^{C_{0}i}\s^{-n-C_{0}i}(1+|v|)^{3i+2n+(i+n-1)\bar\gamma} \mu_{\sigma},\label{3.27-1}
\end{align}
where $\bar{\gamma}=\max{\{0,\gamma\}}$, and $C_0, C_{i}, C_{i,n}\geq1$ depend only on  $\|(\r_{i},u_i,\t_{i})(0)\|_{H^s}$ and $\theta_{\pm}$.
\end{theorem}

\begin{remark}
We can not use the classical results \cite{Caflish} on the linear terms $F_i, i=1,2,3,4,5$ since the $F_0$ depends on the parameter $\sigma$. Indeed, from \eqref{3.26-1} and \eqref{3.27-1}, we know that $F_i$  grows polynomially as $\sigma\rightarrow0+$, and this fact is very important for us to prove the hydrodynamic limit below.
\end{remark}

We shall construct a sequence of solution of  Boltzmann equation near  the local Maxwellian $\mu_{\sigma}$, so it is natural to rewrite the remainder as
\begin{equation}\label{1.20}
F^{\ep}_{R}=\sqrt{\mu_{\sigma}} \fep.
\end{equation}
To use the $L^2$-$L^\infty$ framework \cite{Guo Jang Jiang-1}, we also introduce a global Maxwellian
	\begin{equation*}
		\mu_{M}:=\frac{1}{(2\pi \t_{M})^{3/2}}\exp{\left\{-\frac{|v|^{2}}{2\t_{M}}\right\}},
	\end{equation*}
	where $\t_{M}$ satisfies the condition
	\begin{equation}\label{condition of T}
		\t_{M}<\max_{t\in [0,\infty),x\in\Omega}\t_0(t,x)<2\t_{M}.
	\end{equation}
Since $\theta^{R_1}(t,x_1)$ is a monotonic function of $x_1$, and $\min\{\theta_-,\theta_+\}\leq \theta^{R_1}(t,x_1)\leq \max\{\theta_-,\theta_+\}$, then we can always choose $\theta_{M}$ satisfying \eqref{condition of T} if
\begin{equation}\label{2.16}
\max\{\theta_-,\theta_+\}<2\min\{\theta_-,\theta_+\}.
\end{equation}
By the assumption \eqref{2.16}, one can easily deduce that there exists positive constant $C>0$ such that for some $\frac{1}{2}<\alpha<1$ and for each $(t,x,v)\in [0,\infty)\times \mathbb{R}^3\times \mathbb{R}^{3}$, the following holds:
\begin{equation}\label{relation of mu and muM}
\frac1C\mu_{M}\leq \mu_\sigma(t,x,v)\leq C\mu_{M}^{\alpha},
\end{equation}
where both $C$ and $\alpha$ are independent of $\sigma$.
We further define
\begin{equation}\label{def of h}
F_{R}^{\ep}=\{1+|v|^{2}\}^{-\beta}\sqrt{\mu_{M}}\hep\equiv\frac{1}{w(v)}\sqrt{\mu_{M}}\hep,
\end{equation}
with $w(v):=\{1+|v|^{2}\}^{\beta}$ for any fixed $	\beta\geq \frac{9}{4}+2(3-\gamma)$.  

\begin{theorem}\label{theorem}
Under the assumption of Theorem \ref{thm1.2}, and let \eqref{2.16} hold and $\sigma=\v^{\eta}$, $a=\v^{-2\eta}$. Assume the initial data
\begin{equation}
F^{\ep}(0,x,v)=\mu_\sigma(0,x_1,v)+\sum_{n=1}^{5}\ep^n F_n(0,x_1,v)+\ep^3F^\ep_R(0,x,v)\geq 0,\nonumber
\end{equation}
and
\begin{align}\label{1.36}
\sup_{x_0\in\mathbb{R}^3}\|\fep(0,\cdot,\cdot) I_{\{|\cdot-x_0|\leq 2a\}}\|_{L^2_{x,v}}\lesssim \v^{-\frac18} a^3,\quad \|\frac{\ep^{3/2}}{a^{3}}\hep(0)\|_{L^\infty_{x,v}}\lesssim 1.
\end{align}
Then there are small positive constants $\eta\in(0,\frac1{100})$ and $\ep_0>0$ depending only on $\theta_{\pm}$ such that the Cauchy problem of  Boltzmann equation \eqref{1.1}, \eqref{1.36} has a unique solution for  $\v\in(0,\v_0]$
\begin{equation}\label{1.36-1}
F^{\ep}(t,x,v)=\mu_\sigma(t,x_1,v)+\sum_{n=1}^{5}\ep^n F_n(t,x_1,v)+\ep^3F^\ep_R(t,x,v)\geq 0,\quad t\in[0,\v^{-\eta}],
\end{equation}
with
\begin{align}
\sup_{0\leq t\leq \v^{-\eta}} \sup_{x_0\in\mathbb{R}^3}  \|\fep(t,\cdot,\cdot) I_{\{|\cdot-x_0|\leq 2a\}}\|_{L^2_{x,v}}&\lesssim \ep^{-\frac{33}{200}}a^3,\label{1.38}\\
\sup_{0\leq t\leq \v^{-\eta}} \Big\{\|\frac{\v^{\frac32}}{a^3}h^{\ep}(t)\|_{L^\infty_{x,v}}\Big\}&\lesssim \ep^{-\frac14}.\label{1.37}
\end{align}
\end{theorem}

\begin{remark}
Under the condition \eqref{2.16},  the wave strength of the rarefaction wave may be large in some cases. For example, for 1-rarefaction wave, one can choose $\theta_+=\frac{3}{4}\theta_-$, then it is easy to check that \eqref{2.16} holds. The wave strength $|\theta_+-\theta_-|=\frac14\theta_-$ is large when $\theta_-$ is large.
\end{remark}

\begin{remark}
Since the approximate planar rarefaction wave depends on $\sigma$ (or $\v$), unlike \cite{Guo Jang Jiang}, the uper bound of $L^2$ and $L^\infty$-norms can not be  kept. Indeed, from \eqref{1.37} and \eqref{1.38}, these norm of Boltzmann solution will increase with higher rate than the initial data.
\end{remark}

\begin{remark}
We notice that the functions $\mu_{\sigma}, F_1, \cdots, F_5$ are independent of the space variables $x_2, x_3$. However,   $F^\v(t,x,v)$ is indeed a nontrivial Boltzmann solution in three dimensional space since the remainder term $F^\ep_R(t,x,v)$ depends on $x_1, x_2$ and $x_3$.
\end{remark}

\begin{remark}
Both $\mu_{\sigma}$ and linear terms $F_1,\cdots, F_5$ depend on the $\v>0$ in Theorem \ref{theorem}, hence we call  \eqref{1.36-1} as a generalized Hilbert expansion (The linear part are independent of $\v$ in the classical Hilbert expansion \cite{Caflish,Guo Jang Jiang}).
\end{remark}

\begin{remark}
Under the conditions \eqref{1.36} and \eqref{1.18-2}, one can indeed construct initial data $F_0^\v\geq0$ (we shall not present the details of construction for simplicity), hence the positivity of Boltzmann solution $F^\v(t,x,v)$ can be guaranteed.  
\end{remark}

From \eqref{1.16} and Theorems \ref{thm1.2} and \ref{theorem}, one can obtain the hydrodynamic limit of the nontrivial three dimensional Boltzmann solution to the planar rarefaction wave of compressible Euler equations.
\begin{corollary}[Hydrodynamic limit to the planar rarefaction wave]
Recall the definition of $\mu(t,x_1,v)$ in \eqref{1.18-1}. Under the conditions of Theorem \ref{theorem}, we have the following hydrodynamic limit of Boltzmann equation to the planar rarefaction wave of compressible Euler equations
\begin{align}\nonumber
\sup_{t\in[\v^\zeta,\v^{-\eta}]}\left\|\frac{F^\v(t,x,v)-\mu(t,x_1,v)}{\sqrt{\mu_M}}\right\|_{L^\infty_{x,v}}\lesssim \v^{\eta-\zeta}|\ln\v|\rightarrow0+,\  \mbox{as}\  \v\rightarrow 0+,
\end{align}
for any  given positive constant $\zeta\in(0,\eta)$.
\end{corollary}
\begin{remark}
As pointed out in the introduction, all the results \cite{Yu,Xin-Zeng,Huang-Jiang-Wang,Huang-Wang-Yang,Huang-Wang-Yang-1,Huang-Wang-Wang-Yang} on hydrodynamic limit of Boltzmann equation to the wave pattern solution of Euler system are  one dimensional case, i.e. $x\in\mathbb{R}$. In the present paper, we provide the first result on the hydrodynamic limit of Boltzmann equation to the planar wave pattern solution of compressible Euler system in three dimensional space $x\in\mathbb{R}^3$. On the other hand,
the validity time in the hydrodynamic limit is $\v^{-\eta}$ for some small positive constant $\eta>0$, which implies the global in-time convergence from Boltzmann solution to planar rarefaction wave of the compressible Euler system.
\end{remark}

We now comment on the analysis of this paper. For the linear part $F_1,\cdots, F_5$, we can not use the classical results \cite{Caflish} since the local Maxwellian $\mu_{\sigma}$ (see \eqref{1.18} for definition) depends on the parameter $\sigma>0$, and the linear part $F_i$ may grow to infinity when $\sigma$ vanishes.  Hence one needs to obtain a growth estimation as  $\sigma\rightarrow 0$.  Noting the  properties \eqref{estimate 2.0} of approximate rarefaction wave, one can prove that the linear parts $F_1,\cdots, F_5$ satisfy
\begin{equation*}
|F_i(t,x_1,v)|\leq C_i (\sigma+t)^{C_0i}\s^{-C_0i}(1+|v|)^{3i+(i-1)\bar{\gamma}}\mu_\sigma,
\end{equation*}
which grow to infinity with polynomial rate as $\sigma\rightarrow0+$, see section \ref{section2} for details.

For the estimation of  reminder term $F^\v_{R}$, our method of proof relies on a recent $L^2-L^\infty$ framework initiated in \cite{Guo2010,Guo Jang Jiang-1}. Since the planar rarefaction wave and linear parts are independent of $x_2, x_3$, the source term  $\bar{A}(t,x_1,v)$ defined in \eqref{3.2} is not integrable in $L^2(\mathbb{R}^3_x\times\mathbb{R}^3_v)$. To overcome such difficulty, we introduce a localized $L^{2}_{x,v}$ estimation for $f^\v$. In fact, we consider the equation of $f^\v(t,x,v) \varphi_a(x-x_0)$ for any given $x_0\in\mathbb{R}^3$ to obtain
\begin{align}
&\dis \pt(\fep\varphi_{a})+v\cdot\nax(\fep\varphi_{a})+\frac{1}{\ep}\mathbf{L}(\fep\varphi_{a})\nonumber\\
&=-\frac{\{\pt+v\cdot\nax\}\sqrt{\mu_\sigma}}{\sqrt{\mu_\sigma}}\fep\varphi_{a}+(v\cdot\nabla_x)\varphi_{a}\fep+\ep^{2}\Gamma(\fep,\fep\varphi_{a})\nonumber\\
&\dis  \quad+\sum_{i=1}^{5}\ep^{i-1}\Big\{\Gamma(\frac{F_{i}}{\sqrt{\mu_\sigma}},\fep\varphi_{a})+\Gamma(\fep\varphi_{a},\frac{F_{i}}{\sqrt{\mu_\sigma}})\Big\}+\v^2\bar{A}(t,x_1,v) \varphi_{a},\nonumber
\end{align}
where $\varphi_a$ is a cut-off function on $x$ defined in \eqref{cut-off} and $a=\v^{-2\eta}$. Compared to \cite{Guo Jang,Guo Jang Jiang-1,Guo Jang Jiang}, the term $(v\cdot\nabla_x)\varphi_{a}\fep$ is new. To close the estimate,  we  have to be careful since we need some $\v$ power to macth the term $\|h^\v\|_{L^\infty}$. Noting the definition of $\varphi$ in \eqref{cut-off}, one has 
$$|(v\cdot\nabla_x)\varphi_{a}|\leq C_\lambda a^{-1-3\lambda} |v|\cdot|\varphi_a|^{1-\lambda}\  \mbox{for} \  \lambda\in(0,1),$$
which provides an additional $\v$ decay, i.e.,  $a^{-1}=\v^{2\eta}$. And this is the main reason why we choose the cut-off parameter $a$ to depend on $\v$.
Hence the energy estimate of this term can be bounded as
\begin{align}
\left|\intr\intr v\cdot\nax\varphi_{a}|\fep|^2\varphi_{a} dx dv\right|
&\leq \frac{C_{\lambda}}{a^{1+\frac{3}{2}\lambda}}\|\hep\|_{L^\infty}^{\lambda}\cdot\|\fep\varphi_{a}\|_{L^2}^{2-\lambda} \nonumber\\
&\leq C_\lambda \v^{\frac32\eta}\left(\v^{\frac14}\|\frac{\v^{\frac{3}{2}}}{a^3} h^\v(t)\|_{L^\infty}\right)^{\lambda}\cdot\|\fep\varphi_{a}\|_{L^2}^{2-\lambda},\nonumber
\end{align}
by taking $\lambda=\frac{1}{21}\eta$, see \eqref{4.9} for more details.  We emphasize that the gain of $\v$ power from the $\nabla\varphi_a$ is one of the key point. For the other terms in the energy estimates, one can bound them by similar arguments as in \cite{Guo Jang Jiang-1,Guo Jang Jiang}. Hence, by choosing $\sigma=\v^{\eta}$ with $\eta>0$ being suitably small, we can obtain that
\begin{align}\label{1.41}
&\frac{d}{dt}\|\fep(t)\varphi_{a}(\cdot-x_0)\|_{L^2}^{2}+\frac{c_{0}}{2\ep}\|\{\mathbf{I-P}\}(\fep(t)\varphi_{a}(\cdot-x_0))\|_{\nu}^{2}\nonumber\\
&\leq  \frac{4\tilde{C}_1}{\sigma+t}
\cdot(\|\fep(t)\varphi_{a}(\cdot-x_0)\|_{L^2}^{2}+1),\quad\mbox{for}\  t\in[0,\v^{-\eta}],
\end{align}
where we have used the {\it a priori} assumption $\sup_{0\leq t\leq \v^{-\eta}} \Big\{\v^{\frac14}\|\frac{\v^{\frac{3}{2}}}{a^3} h^\v(t)\|_{L^\infty}\Big\}\leq 1$, see Lemma \ref{lem5.1} and \eqref{4.46}-\eqref{4.48} for details. The key point is that the positive constant $\tilde{C}_1$ is independent of $x_0\in\mathbb{R}^3$. The second step is to estimate the weighted $L^\infty$-norm so that we can close the {\it a priori} assumption, and the key obsevation is that  such local $L^2$-estimate is enough to close the weighted $L^\infty$-estimate, i.e.,
\begin{align}\label{1.42}
&\sup_{0\leq s\leq t}\|\frac{\ep^{3/2}}{a^{3}}\hep(s)\|_{L^\infty}
\leq C\Big\{\|\frac{\ep^{3/2}}{a^{3}}\hep(0)\|_{L^\infty}+C\frac{\ep^{9/2}}{a^{3}}(1+t)^{10C_0}\cdot \sigma^{-10C_0}\Big\}\nonumber\\
&\qquad\qquad\qquad+C\ep^{3/2}a^{3}\sup_{0\leq s\leq t}\|\frac{\ep^{3/2}}{a^{3}}\hep(s)\|^{2}_{L^\infty}+C\sup_{x_0\in\mathbb{R}^3}\sup_{0\leq s\leq t}\|\fep(s)\varphi_{a}(\cdot-x_0)\|_{L^2}.
\end{align}
where have used the fact $t\in[0,\v^{-\eta}]$. With the help of \eqref{1.41}, \eqref{1.42} and the continuity argument, we can finally prove Theorem \ref{theorem}. \vspace{1.5mm}

The paper is organized as follows. In Section 2, we introduce some useful lemmas which will be used later.  In Section 3, we construct the coefficients $F_i$ for the Hilbert expansion for any given $\mu_{\sigma}$, and obtain some estimates depending on  $\sigma$. In Section 4, we derive the localized $L^2$ energy estimate for the remainder $f^\v$ in terms of weighted $L^\infty$-norm, and also  the weighted $L^\infty$-norm in terms of  the localized $L^2$-norm. The main Theorem \ref{theorem} is proved based on the interplay of $L^2$-$L^\infty$ estimates. \vspace{1.5mm}

\noindent{\bf Notations.}  Throughout this paper, $C$ denotes a generic positive constant which may depend on $\rho_{\pm},u_{\pm},\theta_{\pm}$  and  vary from line to line but independent of $\v,\sigma, t$. And $C_a,C_b,\cdots$ denote the generic positive constants depending on $a,~b,\cdots$, respectively, but independent of $\v,\sigma, t$, which also may vary from line to line. $\|\cdot\|_{L^2}$ denotes the standard $L^2(\mathbb{R}^3_x\times\mathbb{R}^3_v)$-norm, and $\|\cdot\|_{L^\infty}$ denotes the $L^\infty(\mathbb{R}^3_x\times\mathbb{R}^3_v)$-norm.\vspace{1.5mm}
	
\section{Preliminaries}
We introduce the following notation
$$\pa_{\tau}^\a=\pa_{t}^{\a_{0}}\pa_{x_1}^{\a_{1}}.$$
We denote $|\a|=\a_{0}+\a_{1}$ where $\a_{0},\a_{1}\in \mathbb{N}, \a_{0},\a_{1}\geq0$. For simplicity, we represent $\pa_{\tau}^\a$ by $\pa^n_\tau$ for the case $|\a|=n$. The properties on the approximate rarefaction wave $(\rho^{R_1}, u^{R_1}, \theta^{R_1})(t,x_1)$ can be summarized as follows.
\begin{lemma}[Xin \cite{Xin-1993}]
The approximate rarefaction waves $(\r^{R_1},u^{R_1}, \theta^{R_1})(t,x_1)$ constructed in \eqref{construct 1 of r u t} have the following properties:\vspace{1mm}
	
	
\noindent (1) For any $1\leq p\leq +\i$ and $k\geq 2$, the following estimates holds,
\begin{equation}\label{estimate 2.0}
	\begin{array}{ll}
	\|\pa_{\tau}(\r^{R_1},u_1^{R_1}, \theta^{R_1})(t,\cdot)\|_{L^p(dx_1)} \leq
	C(\sigma+t)^{-1+\f1p},\quad \\
	\|\pa_{\tau}^\a(\r^{R_1},u_1^{R_1}, \theta^{R_1})(t,\cdot)\|_{L^p(dx_1)} \leq
	C(\sigma+t)^{-1}\cdot \sigma^{-k+1+\f1p},\quad |\a|=k\geq2\\
	\end{array}
\end{equation}
where the positive constant $C$ depends only  on $p,k$ and the wave strength $|\theta_+-\theta_-|$.\vspace{1mm}
	
	
\noindent (2)There exist positive constants $C>0$ and $\s_0>0$ such that for $\s\in(0,\s_0)$ and $t>0,$
\begin{align}
\sup_{x_1\in\mathbb{R}}\left|(\r^{R_1},u_1^{R_1},\theta^{R_1})(t,x_1)-(\r^{r_1},u_1^{r_1},\theta^{r_1})(\f{x_1}t)\right|\leq\f{C}{t}[\s\ln(1+t)+\s|\ln\s|].
\end{align}
\end{lemma}

\

A direct calculation shows that  $\pa_\tau \m_\sigma=\m_{\sigma} J_\tau$ where
\begin{equation}\label{2.6}
J_\tau (t,x_1,v):=\f{\pa_\tau \r_0}{\r_0}-\f32\f{\pa_\tau\t_0}{\t_0}+\f{(v-u_0)\pa_\tau u_0}{\t_0}+\f{|v-u_0|^2 \pa_\tau \t_0}{2\t_0}.
\end{equation}
For $k\geq 0$, it follows from \eqref{estimate 2.0}  that
\begin{align}
|\pa_\tau^k J_\tau|&\leq \sum_{i=0}^{k}C_{k}\s^{-k+i} \cdot| (\pa_\tau^{i+1}\r_0, \pa_\tau^{i+1} u_0, \pa_\tau^{i+}\t_0)|\cdot(1+|v|)^2\nonumber\\
&\leq C_{k}\f{\sigma^{-k}}{\s+t}(1+|v|)^2
\leq C_{k}\s^{-(k+1)}(1+|v|)^2,\label{2.8}
\end{align}
where $C_{k}$ is a constant depending on $k$ and wave strength $|\theta_+-\theta_-|$.

For later use, we introduce some linear spaces, functions and operators. Based on $J_\tau$,
we define the operators $A_k$
\begin{equation}\label{def of A}
A_0(f)=f,\quad A_1(f)=\pa_\tau f+\f12 fJ_\tau,\quad A_{k+1}(f)=A_{1}\circ A_{k}(f),
\end{equation}
and the linear spaces $B_k$
\begin{equation}\label{def of B}
B_1=span\{J_\tau\},\quad B_2=span\{\pa_\tau J_\tau, J_\tau^2\},...,B_k=span\Big\{\prod_{i=1}^{k}\pa_\tau^{m_i-1}J_\tau\Big\}_{|m(k)|=k},
\end{equation}
where $m(k)=(m_1,...,m_k)$ is a multiple index with $m_i\in \mathbb{N}, m_i\geq 0$ and $|m(k)|=k$. 
For later use, we also denote $\pa_\tau^{-1} J_\tau=1, \pa_\tau^0 J_\tau=J_\tau$ and $b_0=1$.

Now we give some useful lemmas which will be used in section \ref{section2}. The proofs of Lemmas \ref{lem4.1}, \ref{prop of A}, \ref{derivatives of \L} and \ref{estimate of Gamma} are presented in the Appendix.
\begin{lemma}\label{lem4.1}
For the linear space $B_i, i\geq 1$, we have the following properties \vspace{1.5mm}\\
\noindent{1)} Let $b_p\in B_p$ and $b_q\in B_q$, then it holds that  $b_p b_q \in B_{p+q}$;\vspace{2mm}

\noindent{2)} Let $b_n\in B_n$ and $f$ be any smooth function, then there exists a $b_{n+1}\in B_{p+1}$ such that
\begin{equation}
\partial_\tau b_p=b_{p+1}\in B_{p+1}\  \mbox{and}\ \  \pa_\tau (b_nf)=b_{n+1}f+b_n\pa_\tau f;
\end{equation}

\noindent{3)} 	Let $b_k\in B_k$ be the basis of $B_k$, i.e. $b_k=\prod_{i=1}^{k}\pa_\tau^{m_i-1}J_\tau$ for some $m(k)$,  then it holds that
\begin{equation}\label{2.12}
|b_k |\leq 
C_{k}\s^{-k}(1+|v|)^{2k},
\end{equation}
where $C_{k}$ is a positive constant depending on $k$ and wave strength $|\theta_+-\theta_-|$.
\end{lemma}


\

Let $f\in \mathcal{N}^{\perp}$ and $b_i\in B_i,  b_j\in B_j$, we define a new operators
\begin{equation}\label{2.9}
\Gamma_{i,j}(f)=\frac{1}{\sqrt{\mu_{\s}}}\Big[Q(b_i\mu_{\s}^{\f12}\mathbf{L}^{-1}f,b_j\mu_{\s})+Q(b_j\mu_{\s}, b_i\mu_{\s}^{\f12}\mathbf{L}^{-1}f)\Big],\quad i\geq 0, j\geq 1,
\end{equation}
and we also define $\Gamma_{0,0}=id$. For simplicity of presentation, we may still use the same notation $\Gamma_{i,j}(f)$ even though $b_i\in B_i, b_j\in B_j$ are replaced by other $\tilde{b}_i\in B_i, \tilde{b}_j\in B_j$. And such simplification will not cause problem in the following estimations.
\begin{lemma}\label{prop of A}
1). There exist $b_i\in B_i, i=0,\cdots, k$ such that 
\begin{align}\label{prop of A_k}
A_kf&=\sum_{i=0}^kb_i\pa_\tau^{k-i}f,\quad k\geq1.
\end{align}

2). There exist $b_1\in B_1, b_{i+1}\in B_{i+1}$ and $b_{j+1}\in B_{j+1}$ such that 
\begin{align}
A_1\circ \Gamma_{i,j}(f)&=b_1\Gamma_{i,j}(f)+\Gamma_{i+1,j}(f)+\Gamma_{i,j+1}(f)+\Gamma_{i,j}\circ \Gamma_{0,1}(f)\nonumber\\
&\quad+\Gamma_{i,j}\circ A_1(f)\quad i\geq 0, j\geq 1.\label{A and G}
\end{align}
\end{lemma}

\begin{lemma}\label{derivatives of \L}
There exist index sets $$N_s(i,j,l)=\Big\{(i_m,j_m,l_m)\in \mathbb{N}_+^3|\sum_{m=1}^{s}(i_m,j_m,l_m)=(i,j,l), \mbox{and} \   i_m=j_m=0\  \mbox{as}  \  l_m=0\Big\},$$
such that
\begin{equation}\label{n-th derivates of L^-1}
\pa_\tau^n \mathbf{L}^{-1} f =\sum_{\substack{r+k=n\\r,k\geq0}}\sum_{\substack{s+p=k\\s,p\geq0}}b_r\mathbf{L}^{-1}[\sum_{\substack{i+j+l=s\\i,j,l\geq0}}\sum_{(i_m,j_m,l_m)\in N_s(i,j,l)}\left(b_{i_1}\Gamma_{j_1,l_1}\right)\circ\cdots\circ\left(b_{i_s}\Gamma_{j_s,l_s}\right)\circ A_pf],
\end{equation}
for any $f\in \mathcal{N}^{\bot}$, $n\geq0$.
\end{lemma}

\begin{lemma}\label{estimate of Gamma}
Let $f\in\mathcal{N}^{\perp}$ and $|f(t,x,v)|\leq S(t,x)(1+|v|)^m\sqrt{\mu_{\s}}$ where $S(t,x)\geq 0$, then it holds that
\begin{equation}\label{2.14}
|\Gamma_{i,j}(f)|\leq C_{i,j}\s^{-(i+j)}S(x,t)(1+|v|)^{m+2i+2j+\gamma}\sqrt{\mu_{\s}},
\end{equation}
where $C_{i,j}$ is positive constant independent of $\sigma$.
\end{lemma}

\vspace{1.5mm}

\section{Estimates on the linear terms}\label{section2}

In this section, we will derive the estimates of $F_1(t,x_1,v),\cdots F_5(t,x_1,v)$ for given $\mu_{\sigma}$ which is defined \eqref{1.18}.
We also point out that all the functions are independent of $x_2$ and $x_3$ in this section.
We define $\displaystyle f_k:=\frac{F_k}{\sqrt{\mu_{\sigma}}}$. Firstly, we present a useful lemma in \cite{Guo Jang} which will be used to estimate the bound of linear terms.
\begin{lemma}[Guo-Jang \cite{Guo Jang}]\label{lem2.1}
	For each given nonnegative integer $k$, assume $f_k$'s are found. Then the microscopic part of $f_{k+1}$ is determined through the equation for $F_k$ in  \eqref{F_k equations}:
	\begin{equation}\label{(I-P)f_k}
	\dis \{\mathbf{I-P}\}f_{k+1}=\mathbf{L}^{-1}\left(-\f{\{\pa_t+v^1\pa_{x_1}\}(\sqrt{\mu_{\s}} f_k)-\sum_{\substack{i+j=k+1\\i,j\geq 1}}Q(\sqrt{\mu_{\s}} f_i,\sqrt{\mu_{\s}} f_j)}{\sqrt{\mu_{\s}}}\right).
	\end{equation}
	For the macroscopic part, $\r_{k+1},u_{k+1},\t_{k+1}$ satisfy the following:
	\begin{equation}\label{macro equations}
	\begin{array}{c}
	\dis \pa_t\r_{k+1}+\pa_{x_1}(\r_0u_{k+1}^1+\r_{k+1}u_0^1)=0,\\[1.5mm]
	\dis \r_0\Big\{\pa_t u_{k+1}^{1}+u_{k+1}^{1}\pa_{x_1} u_0^1+u_0^1\pa_{x_1}u_{k+1}^1\Big\}-\f{\r_{k+1}}{\r_0}\pa_{x_1}(\r_0\t_0)+\pa_{x_1}(\f{\r_0\t_{k+1}+3\t_0\r_{k+1}}{3})=\bar{f}_{k,1},\\[3mm]
	\dis \r_0\{\pa_t u_{k+1}^{2}+u_0^1\pa_{x_1}u_{k+1}^2\}=\bar{f}_{k,2},\\[1.5mm]
	\dis \r_0\{\pa_t u_{k+1}^{3}+u_0^1\pa_{x_1}u_{k+1}^3\}=\bar{f}_{k,3},\\[1.5mm]
	\dis \rho_0\Big\{\pa_t\t_{k+1}+\f23(\t_{k+1}\pa_{x_1}u_0^1+3\t_0\pa_{x_1}u_{k+1}^{1})+u_0^1\pa_{x_1}\t_{k+1}+3u_{k+1}^1\pa_{x_1}\t_0\Big\}=\bar{g}_k,
	\end{array}
	\end{equation}
where
	\begin{align}
	\begin{split}
	\bar{f}_{k,i}&=-\pa_{x_1}\left(\t_0\int_{\mathbb{R}^3}\mathcal{B}_{i,1}F_{k+1}dv\right),\\
	\bar{g}_{k}&=-\pa_{x_1}\left(\t_0^{\f32}\int_{\mathbb{R}^3}\mathcal{A}_{1}F_{k+1}dv+2u_0^1\t_0\int_{\mathbb{R}^3}\mathcal{B}_{1,1}F_{k+1}dv\right)-2u_0^1f_{k,1},
	\end{split}
	\end{align}
and
	\begin{equation}\nonumber
	\mathcal{A}_{i}=\f{v^i-u^i_0}{\sqrt{\t_0}}\left(\f{|v-u_0|^2}{\t_0}-5\right),\quad \mathcal{B}_{i,j}=\f{(v^i-u_0^i)(v^j-u_0^j)}{\t_0}-\d_{ij}\f{|v-u_0|^2}{3\t_0},
	\end{equation}
where we use the subscript $k$ for forcing terms $\bar{f}_{k,i}$ and $\bar{g}_{k,i}$ in order to emphasize that the right hand side depends only on $F_i$'s for $0\leq i\leq k$.
\end{lemma}

\begin{remark}
	The original version of Lemma \ref{lem2.1} in \cite{Guo Jang} is for the Hilbert expansion of Vlasov-Poisson-Boltzmann equations,  and one can obtain Lemma \ref{lem2.1} by dropping the electric field and noting that all the functions are independent of variables $x_2$ and $x_3$.
\end{remark}


\noindent{\bf Proof of Theorem \ref{thm1.2}.} 	Firstly we consider the microscopic part $\{\mathbf{I-P}\}f_1$. It follows from \eqref{(I-P)f_k} that
\begin{equation}\label{6.4} \{\mathbf{I-P}\}f_1=\mathbf{L}^{-1}\left(-\f{\{\pa_{t}+v^{1}\pa_{x_1}\}\mu_\s}{\sqrt{\mu_\s}}\right)=\mathbf{L}^{-1}(- J_t\sqrt{\mu_\s}-v^1J_{x_1}\sqrt{\mu_{\s}}).
\end{equation}
Since $\mathbf{L}^{-1}$ preserves decay of $v$ \cite{Caflish}, and $\r_0, \t_0$ are  bounded from  below and above, then using \eqref{2.8} to obtain
\begin{align}\label{6.3}
|\{\mathbf{I-P}\}f_1|& \leq C|(\pa_\tau\r_0, \pa_\tau u_0, \pa_\tau\t_0)|\cdot(1+|v|)^3\sqrt{\mu_{\s}}\nonumber\\
&\leq C\sigma^{-1}(1+|v|)^3\sqrt{\mu_{\s}}.
\end{align}
It follows from  \eqref{6.3} and \eqref{estimate 2.0} that
\begin{align}\label{6.7}
&\|\{\mathbf{I-P}\}f_1\|_{L^2_{x_1}L^\i_v}+\|\{\mathbf{I-P}\}f_1\|_{L^2_{x_1}L^2_v}\nonumber\\
&\leq C\|(\pa_{\tau}\r_0, \pa_{\tau}u_0, \pa_{\tau}\t_0)\|_{L^2_{x_1}}\leq C(\sigma+t)^{-\frac12}\leq C \s^{-\f12}.
\end{align}

Next we consider the space-time derivatives of $\{\mathbf{I-P}\}f_1$. It follows from \eqref{6.4} and Lemma \ref{derivatives of \L} that
\begin{align}\label{6.8}
\pa_\tau^n\{\mathbf{I-P}\}f_1&=\pa_\tau^n\mathbf{L}^{-1}(- J_t\sqrt{\mu_{\s}}-v^1J_{x_1}\sqrt{\mu_{\s}})\nonumber\\
&=\sum_{\substack{r+k=n\\r,k\geq0}}\sum_{\substack{s+p=k\\s,p\geq0}}b_r\mathbf{L}^{-1}\Big[\sum_{\substack{i+j+l=s\\i,j,l\geq0}}\sum_{(i_m,j_m,l_m)\in N_s(i,j,l)}\left(b_{i_1}\Gamma_{j_1,l_1}\right)\circ\nonumber\\
&\qquad\qquad\qquad\cdots\circ\left(b_{i_s}\Gamma_{j_s,l_s}\right)\circ A_p(- J_t\sqrt{\mu_\s}-v^1J_x\sqrt{\mu_{\s}})\Big].
\end{align}
It is noted that there exists some $b_k\in B_k$ (see \eqref{def of B} for the definition of $B_k$) such that $\pa_\tau^k\sqrt{\mu_{\s}}=b_k\sqrt{\mu_{\s}}$ for $\ k\geq 0$. By using Lemmas \ref{lem4.1} and \ref{prop of A}, it holds that
\begin{align}\label{6.9}
&A_p( J_t\sqrt{\mu_{\s}}+v^1J_{x_1}\sqrt{\mu_{\s}})\nonumber\\
&=\sum_{i=0}^p\sum_{j=0}^{p-i}b_i\pa_\tau^{p-i-j}\sqrt{\mu_\s} \cdot\pa_\tau^jJ_t+v^1\sum_{i=0}^p\sum_{j=0}^{p-i}b_i\pa_\tau^{p-i-j}\sqrt{\mu_{\s}}\cdot\pa_\tau^jJ_{x_1}\nonumber\\
&=\sum_{j=0}^{p} \Big[b_{p-j}\pa_\tau^jJ_t+v^1b_{p-j}\pa_\tau^jJ_{x_1}\Big]\sqrt{\mu_{\s}},\nonumber\\
&=(b_{p+1}+v^1\tilde{b}_{p+1}) \sqrt{\mu_{\s}}.
\end{align}
Substituting \eqref{6.9} into \eqref{6.8}, then using \eqref{estimate 2.0} and \eqref{2.12},  Lemmas \ref{derivatives of \L} and \ref{estimate of Gamma}, one obtains that
\begin{align}
\left|\pa_\tau^n\{\mathbf{I-P}\}f_1\right|
&\leq C_n\sum_{\bar{j}=0}^{n}\s^{-n+\bar{j}}(1+|v|)^{2n+3+n\bar{\gamma}}\cdot|\pa_\tau^{\bar{j}+1}(\r_0, u_0, \t_0)|\sqrt{\mu_{\s}}\label{6.10}\\
&\leq C_n\s^{-n-1}(1+|v|)^{2n+3+n\bar\gamma}\sqrt{\mu_{\s}},\  \mbox{for}\ n\geq 1, \label{6.11}
\end{align}
where $\bar{\gamma}=\max{\{\gamma,0\}}$.  Using \eqref{6.10} and \eqref{estimate 2.0}, it holds that for $n\geq1$
\begin{align} \label{6.12}
\|\pa_\tau^n\{\mathbf{I-P}\}f_1\|_{L^2_{x_1}\liv}+\|\pa_\tau^n\{\mathbf{I-P}\}f_1\|_{L^2_{x_1}\ltv}
\leq  C\s^{-n-\f12}.
\end{align}

To estimate $\mathbf{P}f_1$, we rewrite the linear system \eqref{macro equations} as a symmetric hyperbolic equations with the corresponding symmetrizer $\bar{A}_0$
\begin{equation}\label{macro equations 2}
\bar{A}_0\pa_t U_{k+1}+\bar{A}_1\pa_{x_1} U_{k+1}+\bar{B}U_{k+1}=\bar{F}_k,
\end{equation}
where $U_{k+1}=(\r_{k+1}, u_{k+1}, \t_{k+1})^t$, and $\bar{A}_0$, $\bar{A}_1$, $\bar{B}$ and $\bar{F}_k$ are given by
\begin{equation}\nonumber
\quad \bar{A}_0\equiv \left(
\begin{array}{ccccc}
(\t_0)^2 & 0 & 0 & 0 & 0 \\
0 & (\r_0)^2\t_0 & 0 & 0 & 0\\
0 & 0 & (\r_0)^2\t_0 & 0 & 0 \\
0 & 0 & 0 & (\r_0)^2\t_0 & 0 \\
0 & 0 & 0 & 0 & \f{(\r_0)^2}{6} \\
\end{array}
\right),
\end{equation}
\begin{equation}\nonumber
\bar{A}_1\equiv\left(
\begin{array}{ccccc}
(\t_0)^2u_0^1 & \r_0(\t_0)^2 & 0 & 0 & 0 \\
\r_0(\t_0)^2 & (\r_0)^2\t_0u^1_0 & 0 & 0 & \f{(\r_0)^2\t_0}{3} \\
0 & 0 & (\r_0)^2\t_0u^1_0 & 0 & 0 \\
0 & 0 & 0 & (\r_0)^2\t_0u^1_0 & 0 \\
0 & \f{(\r_0)^2\t_0}{3} & 0 & 0 & \f{(\r_0)^2u_0^1}{6} \\
\end{array}
\right),
\end{equation}
and
\begin{equation}\nonumber
\quad \bar{B}\equiv\left(
\begin{array}{ccccc}
(\t_0)^2\pa_{x_1}u_0^1 & (\t_0)^2\pa_{x_1}\r_0 & 0\quad  & 0 & 0 \\
\t_0\pa_{x_1}(\r_0\t_0) & \r_0\t_0\pa_{x_1}u^1_0 & 0\quad  & 0 & \f{\r_0\t_0\pa_{x_1}\r_0}{3} \\
0 & 0 & 0 \quad & 0 & 0 \\
0 & 0 & 0 \quad & 0 & 0 \\
0 & \f{(\r_0)^2\partial_{x_1}\t_0}{2} & 0 \quad & 0 & \f{(\r_0)^2\pa_{x_1}u_0^1}{9} \\
\end{array}
\right),
\quad \bar{F}_k\equiv\left(
\begin{array}{c}
0 \\
\r_0\t_0\bar{f}_{k,1} \\
\r_0\t_0\bar{f}_{k,2} \\
\r_0\t_0\bar{f}_{k,3} \\
\f{(\r_0)^2}{6}\bar{g}_k \\
\end{array}
\right).
\end{equation}

Using \eqref{estimate 2.0}, it is easy to know that
\begin{equation}\label{A_1 and B}
\|\pa_{x_1}(\bar{A}_0,\bar{A}_1)\|_{L^\i_{x_1}}+\|\bar{B}\|_{L^\i_{x_1}}\leq \f{\bar{C}}{\s+t},
\end{equation}
where $\bar{C}\geq 1$ is a positive constant depending only on $\theta_{\pm}$.
Applying the standard energy method of the linear symmetric hyperbolic system to \eqref{macro equations 2} and using \eqref{A_1 and B}, then one can obtain the following energy inequality
\begin{align}\label{k+1-th energy estimate}
\dis \f{d}{dt}\|U_{k+1}\|^2_{L^2_{x_1}}
&\leq  \Big\{\big[\|\pa_{x_1}(\bar{A}_0,\bar{A}_1)\|_{L^\i_{x_1}}+\|\bar{B}\|_{L^\i_{x_1}}\big]\|U_{k+1}\|^{2}_{L^2_{x_1}}+\|\bar{F}_k\|_{L^2_{x_1}}\|U_{k+1}\|_{L^2_{x_1}}\Big\}\nonumber\\
&\leq \f{\bar{C}}{\s+t} \|U_{k+1}\|^{2}_{L^2_{x_1}}+C\|\bar{F}_k\|_{L^2_{x_1}}\|U_{k+1}\|_{L^2_{x_1}}.
\end{align}
To estimate $\|\bar{F}_k\|_{L^2_{x_1}}$, we only calculate the term $\|\r_0\t_0\bar{f}_{k,i}\|_{L^2_{x_1}}$ since all the other terms  can be bounded in a similar way.
Noting that $\dis\int_{\mathbb{R}^3}\mathcal{B}_{i,j}\sqrt{\mu_{\s}}\cdot \mathbf{P}f_1 dv=0$ (see \cite{Bardos} for more details), one has that
\begin{align}
\r_0\t_0\bar{f}_{k,i}
&=-\r_0\t_0\int_{\mathbb{R}^3}\pa_{x_1}(\t_0\mathcal{B}_{i,1}\sqrt{\mu_{\s}})\cdot\{\mathbf{I-P}\}f_{k+1}dv\nonumber\\ &\quad-\r_0(\t_0)^2\int_{\mathbb{R}^3}\mathcal{B}_{i,1}\sqrt{\mu_{\s}}\cdot\pa_{x_1}(\{\mathbf{I-P}\}f_{k+1})dv, \nonumber
\end{align}
which, together with  \eqref{estimate 2.0}, yields that
\begin{align}
\|\r_0\t_0\bar{f}_{k,i}\|_{L^2_{x_1}}
&\leq  C\|\pa_{x_1}(\t_0\mathcal{B}_{i,1}\sqrt{\mu_{\s}})\|_{L^\i_{x_1}\ltv}\cdot\|\{\mathbf{I-P}\}f_{k+1}\|_{L^2_{x_1}\ltv}\nonumber\\
&\quad + C\|\mathcal{B}_{i,1}\sqrt{\mu_{\s}}\|_{L^\i_{x_1}\ltv}\cdot\|\pa_{x_1}\{\mathbf{I-P}\}f_{k+1}\|_{L^2_{x_1}\ltv}.\nonumber
\end{align}
It follows from \eqref{estimate 2.0} that
\begin{equation}\nonumber
\dis \|\mathcal{B}_{i,1}\sqrt{\mu_{\s}}\|_{L^\i_{x_1}\ltv}\leq C, \quad  \|\pa_{x_1}(\t_0\mathcal{B}_{i,1}\sqrt{\mu_{\s}})\|_{L^\i_{x_1}\ltv}\leq \f{C}{\s+t},
\end{equation}
which yields immediately that
\begin{equation}
\|\r_0\t_0\bar{f}_k\|_{L^2_{x_1}}\leq \f{C}{\s+t}\|\{\mathbf{I-P}\}f_{k+1}\|_{L^2_{x_1}\ltv}+C\|\pa_{x_1}\{\mathbf{I-P}\}f_{k+1}\|_{L^2_{x_1}\ltv}.\nonumber
\end{equation}
Hence, by similar arguments, one can prove that
\begin{equation}\label{F}
\dis \|\bar{F}_k\|_{L^2_{x_1}}\leq \f{C}{\s+t}\|\{\mathbf{I-P}\}f_{k+1}\|_{L^2_{x_1}\ltv}+C\|\pa_{x_1}\{\mathbf{I-P}\}f_{k+1}\|_{L^2_{x_1}\ltv}.
\end{equation}
For $k=0$, substituting \eqref{6.7} and \eqref{6.12} into \eqref{F} to have
\begin{equation}\nonumber
\dis \|\bar{F}_0\|_{L^2_{x_1}}\leq 
C\s^{-\f32}
\end{equation}
which, together with \eqref{k+1-th energy estimate}, yields that
\begin{equation}\label{6.20}
\dis \f{d}{dt}\|U_{1}\|^2_{L^2_{x_1}}\leq \f{2\bar{C}}{\s+t}\|U_{1}\|^{2}_{L^2_{x_1}}+C(\sigma+t)\s^{-3},
\end{equation}
where $\bar{C}\geq 1$ is some positive constant which depends only on the wave strength.
Applying the Gronwall's inequality to \eqref{6.20}, then one obtains that
\begin{equation}\label{6.21}
\|U_{1}\|^2_{L^2_{x_1}}(t)\leq C\left(\f{\s+t}{\s}\right)^{2\bar{C}}(U_1^2(0)+\s^{-1})\leq C\left(\f{\s+t}{\s}\right)^{2\bar{C}}\s^{-1},
\end{equation}
where we have used the fact
\begin{align}\label{6.21-1}
\int_0^t (\sigma+\tau)\cdot \left(\frac{\sigma}{\sigma+\tau}\right)^{2\bar{C}}d\tau\leq \frac{\sigma^2}{2\bar{C}-2},\  \mbox{for}\ \  \bar{C}>1.
\end{align}

Next we shall estimate the derivatives of $U_1$. We introduce the following notation
\begin{equation}\nonumber \|\nabla^n\cdot\|_{L^2_{x_1}}=\sum_{\substack{\a_0+\a_1=n\\\a_0,\a_1\geq0}}\|\pa_t^{\a_0}\pa_{x_1}^{\a_1}\cdot\|_{L^2_{x_1}},
\end{equation}
for simplicity of presentation.
Applying $\pa_\tau^\alpha$ to \eqref{macro equations 2} for $k=0$, using the standard energy method to the resultant equation and adding them together for $|\alpha|=n$, then we obtain
\begin{align}\label{6.22}
\f{d}{dt}\|\nabla^n U_1\|_{L^2_{x_1}}^2
&\leq \frac{\bar{C}}{\s+t} \|\nabla^nU_1\|_{L^2_{x_1}}^2+C\sum_{i=2}^n\|\nabla^i( \bar{A}_0,\bar{A}_1)\|_{L^\i_{x_1}}\|\nabla^{n-i+1}U_1\|_{L^2_{x_1}}\|\nabla^nU_1\|_{L^2_{x_1}}\nonumber\\
&\quad +C\sum_{i=1}^n\|\nabla^i \bar{B}\|_{L^\i_{x_1}}\|\nabla^{n-i}U_1\|_{L^2_{x_1}}\|\nabla^nU_1\|_{L^2_{x_1}}\nonumber\\
&\quad +C\|\nabla^n\bar{F}_0\|_{L^2_{x_1}}\|\nabla^nU_1\|_{L^2_{x_1}}.
\end{align}
By using \eqref{estimate 2.0}, a direct calculation shows that
\begin{equation}\label{6.23}
\begin{split}
\|\nabla^i (\bar{A}_0,\bar{A}_1)\|_{L^\i_{x_1}}&\leq \f{C}{\s+t}\s^{-i+1},\ \mbox{for}\ i\geq1, \\
 \|\nabla^i \bar{B}\|_{L^\i_{x_1}}&\leq\f{C}{\s+t}\s^{-i},\ \mbox{for}\ i\geq0.
\end{split}
\end{equation}
For the estimate of $\nabla^n\bar{F}_0$, we only consider the effect of  $\nabla^n(\r_0\t_0\bar{f}_{k,i})$ since the other terms can be done by similar way. In fact, it follows from \eqref{estimate 2.0}, \eqref{6.7} and \eqref{6.12} that
\begin{align}
&\|\nabla^n\bar{F}_0\|_{L^2_{x_1}}\nonumber\\
&\leq \sum_{0\leq i+j\leq n}\|\nabla^i(\r_0\t_0)\|_{L^\i_{x_1}}\Big\{\|\nabla^{1+j}(\t_0\mathcal{B}_{1,1}\sqrt{\m_{\s}})\|_{L^\i_{x_1}\ltv}\|\nabla^{n-i-j}\{\mathbf{I-P}\}f_1\|_{L^2_{x_1}\ltv}\nonumber\\
&\quad +\|\nabla^{j}(\t_0\mathcal{B}_{1,1}\sqrt{\m_{\s}})\|_{L^\i_{x_1}\ltv}\|\nabla^{1+n-i-j}\{\mathbf{I-P}\}f_1\|_{L^2_{x_1}\ltv}\Big\}\nonumber\\
&\leq 
C\s^{-n-\f32},\nonumber
\end{align}
which, together with  \eqref{6.22}, \eqref{6.23} and Cauchy inequality, yields that
\begin{align}\label{6.24}
\f{d}{dt}\|\nabla^n U_1\|_{L^2_{x_1}}^2
&\leq \f{2\bar{C}}{\s+t}\|\nabla^n U_1\|_{L^2_{x_1}}^2+\sum_{i=1}^n C\f{\s^{-2i}}{\s+t}\|\nabla^{n-i}U_1\|^2_{L^2_{x_1}}\nonumber\\
&\quad + C(\sigma+t) \sigma^{-2n-3}. 
\end{align}
For $n=1$, it follows from \eqref{6.24} and \eqref{6.21} that
\begin{align}\nonumber
\f{d}{dt}\|\nabla U_1\|_{L^2_{x_1}}^2
\leq \f{2\bar{C}}{\s+t}\|\nabla U_1\|_{L^2_{x_1}}^2+ \f{C\s^{-3}}{\s+t}\left(\frac{\sigma+t}{\sigma}\right)^{2\bar{C}}
+C(\s+t) \s^{-5},
\end{align}
which, together with Gronwall's inequality and \eqref{6.21-1}, yields that
\begin{align}\label{6.25}
\|\nabla U_1(t)\|^2_{L^2_{x_1}}\leq C \sigma^{-3} \left(\frac{\sigma+t}{\sigma}\right)^{2\bar{C}+\frac12}.
\end{align}

We shall use induction argument to prove that
\begin{align}\label{6.26}
\|\nabla^n U_1(t)\|^2_{L^2_{x_1}}\leq C \sigma^{-2n-1} \left(\frac{\sigma+t}{\sigma}\right)^{2\bar{C}+\frac12},\  \mbox{for}\  n\geq0.
\end{align}
In fact, for $n=0,1$, \eqref{6.26} has already been proved in \eqref{6.21} and \eqref{6.25}. Now we assume that \eqref{6.26} holds for $n\leq k-1$. We consider the case for $n=k$, and it follows from \eqref{6.24} and \eqref{6.26} for $n=1,\cdots, k-1$ that
\begin{align}
\f{d}{dt}\|\nabla^k U_1\|_{L^2_{x_1}}^2
&\leq \f{2\bar{C}}{\s+t}\|\nabla^k U_1\|_{L^2_{x_1}}^2+\sum_{i=1}^k C\f{\s^{-2i}}{\s+t} \sigma^{-2(k-i)-1} \left(\frac{\sigma+t}{\sigma}\right)^{2\bar{C}+\frac12} \nonumber\\
&\quad + C\s^{-2k-3}(\s+t)\nonumber\\
&\leq \f{2\bar{C}}{\s+t}\|\nabla^k U_1\|_{L^2_{x_1}}^2+ C\f{\s^{-2k-1}}{\s+t}\left(\frac{\sigma+t}{\sigma}\right)^{2\bar{C}+\frac12}+C\s^{-2k-3}(\s+t),\nonumber
\end{align}
which, together with Gronwall's inequality, yields
\begin{align}
\|\nabla^k U_1(t)\|_{L^2_{x_1}}^2&\leq C\left\{\|\nabla^k U_1(0)\|_{L^2_{x_1}}^2 + \sigma^{-2k-1} (\frac{\s+t}{\s})^{\frac12}+\s^{-2k-1} \right\}\left(\frac{\sigma+t}{\sigma}\right)^{2\bar{C}}\nonumber\\
\label{7.0.1}&\leq C\sigma^{-2k-1}\left(\frac{\sigma+t}{\sigma}\right)^{2\bar{C}+\frac12}.
\end{align}
Thus we proved \eqref{6.26} holds for $n=k$. Hence \eqref{6.26} holds for $n\geq 0$.

It follows from \eqref{6.26} and Sobolev inequality that
\begin{align}
|\partial^n U_1(t,x_1)|
&\lesssim  \{\|\partial^n U_1(t)\|^2_{L^2_{x_1}}\cdot \|\partial^n \partial_{x_1}U_1(t)\|^2_{L^2_{x_1}}\}^{\frac14}\nonumber\\
&\lesssim \sigma^{-n-1} \left(\frac{\sigma+t}{\sigma}\right)^{\bar{C}+\frac14},\nonumber
\end{align}
which, together with \eqref{1.25}, \eqref{6.11}, yields \eqref{3.26-1} and \eqref{3.27-1} for $i=1$ by suitably chosen $C_0\geq1$.

One can prove \eqref{3.26-1} and \eqref{3.27-1} for $F_2, \cdots, F_5$ step by step by using  similar arguments as for $F_1$ previously, and we omit the details for simplicity of presentation. Therefore the proof of Theorem \ref{thm1.2} is completed. $\hfill\Box$

\section{Proof of the main theorem}
\subsection{Localized $L^2$-estimate} Recalling the definition of $f^\v$ in \eqref{1.20},  we can rewrite the equation \eqref{Remainder of FEP} in terms of $\fep$  as
\begin{align}\label{Rof}
\dis &\pt\fep+v\cdot\nax\fep+\frac{1}{\ep}\mathbf{L}\fep\nonumber\\
&=-\frac{\{\pt+v\cdot\nax\}\sqrt{\mu_\sigma}}{\sqrt{\mu_\sigma}}\fep+\ep^{2}\Gamma(\fep,\fep)\nonumber\\
\dis  &\quad+\sum_{i=1}^{5}\ep^{i-1}\Big\{\Gamma(\frac{F_{i}}{\sqrt{\mu_\sigma}},\fep)+\Gamma(\fep,\frac{F_{i}}{\sqrt{\mu_\sigma}})\Big\}+\v^2 \bar{A}(t,x_1,v),
\end{align}
where
\begin{align}\label{3.2}
 \bar{A}(t,x_1,v)=\sum_{\substack{i+j\geq 6\\1\leq i,j\leq 5}}\ep^{i+j-6}\frac{1}{\sqrt{\mu_{\sigma}}}Q(F_{i},F_{j})-\frac{\{\partial_t+v_1 \partial_{x_1}\} F_5}{\sqrt{\mu_{\sigma}}}.
\end{align}
The last term $\bar{A}(t,x_1,v)$ in \eqref{Rof} is only functions of $x_1$, and it is not integrable in $\mathbb{R}^3$. The key observation is that  only a local $L^2$-estimate is involved when we consider the $L^\infty$ estimation. So to overcome the difficulty, we  consider a localized $L^2$ estimate for $f^\v$.  For later use, we introduce a cut-off function
\begin{equation}\label{cut-off}
\varphi(x)=
\begin{cases}
e^{\f{1}{|x|^2-1}},\quad &|x|<1,\\[1mm]
\quad 0,&|x|\geq 1,
\end{cases}
\end{equation}
and denote $\varphi_a(x)=a^{-3}\varphi(\f x a)$.

\begin{lemma}\label{lem5.1}
Let  $C_0$ be the positive constant defined in Theorem \ref{thm1.2}. Let $\beta\geq \frac94+2(3-\gamma)$,  $\sigma=\v^{\eta}$ with $\eta\leq \frac{1}{11C_0}$ and $T\leq \v^{-\frac{1}{10C_0}}$. Then there exists a suitably small constant $\ep_0>0$  such that for all $\ep\in (0,\ep_0)$, and any fixed $x_0\in\mathbb{R}^3$, it holds that
\begin{align}\label{f L2}
&\quad \frac{d}{dt}\|\fep(t)\varphi_{a}(\cdot-x_0)\|_{L^2}^{2}+\frac{c_{0}}{2\ep}\|\{\mathbf{I-P}\}(\fep(t)\varphi_{a}(\cdot-x_0))\|_{\nu}^{2}\nonumber\\
&\leq \Big\{ \tilde{C}_1 a^2\ep^{2}\sigma^{-\f12}\|\hep(t)\|_{L^\infty}+\frac{C_\lambda}{a^{1+\frac{3}{2}\lambda}}\|\hep(t)\|_{L^\infty}^\lambda\nonumber\\
&\qquad\quad+\tilde{C}_1(1+t)^{10C_0} \v \sigma^{-10C_0}+\frac{\tilde{C}_1}{\sigma+t}\Big\}
\cdot(\|\fep(t)\varphi_{a}(\cdot-x_0)\|_{L^2}^{2}+1),
\end{align}
for $t\in[0,T]$, where $\tilde{C}_1\geq1$ is  positive constant, and   $\lambda>0$ is some small parameter chosen later.
\end{lemma}

\noindent{\bf Proof.} For simplicity of presentation, we only consider the case $x_0=0\in \mathbb{R}^3$ since the proof is the same for $x_0\neq0$. Multiplying \eqref{Rof} by the cut-off function $\varphi_a$, one obtains that
\begin{align}\label{4.5}
&\dis \pt(\fep\varphi_{a})+v\cdot\nax(\fep\varphi_{a})+\frac{1}{\ep}\mathbf{L}(\fep\varphi_{a})\nonumber\\
&=-\frac{\{\pt+v\cdot\nax\}\sqrt{\mu_\sigma}}{\sqrt{\mu_\sigma}}\fep\varphi_{a}+(v\cdot\nabla_x)\varphi_{a}\fep+\ep^{2}\Gamma(\fep,\fep\varphi_{a})\nonumber\\
&\dis  \quad+\sum_{i=1}^{5}\ep^{i-1}\Big\{\Gamma(\frac{F_{i}}{\sqrt{\mu_\sigma}},\fep\varphi_{a})+\Gamma(\fep\varphi_{a},\frac{F_{i}}{\sqrt{\mu_\sigma}})\Big\}+\v^2\bar{A}(t,x_1,v) \varphi_{a}.
\end{align}
Then we multiply \eqref{4.5} by $f^\v\varphi_{a}$ to obtain that
\begin{align}\label{4.3}
&\frac12 \frac{d}{dt} \|f^{\v} \varphi_{a}\|_{L^2}^2+\frac{c_0}{\v} \|\{\mathbf{I-P}\}(f^{\v}\varphi_a)\|_{\nu}^2\nonumber\\
&\leq -\intr\intr \frac{\{\pt+(v\cdot\nax)\}\sqrt{\mu_\sigma}}{\sqrt{\mu_\sigma}} |\fep\varphi_{a}|^2 dv dx+\intr\intr (v\cdot\nax)\varphi_{a} |\fep|^2\varphi_{a}dvdx\nonumber\\
&\quad +\sum_{i=1}^{5}\intr\intr \v^{i-1}\Big\{\Gamma(\frac{F_{i}}{\sqrt{\mu_\sigma}},\fep\varphi_{a})+\Gamma(\fep\varphi_{a},\frac{F_{i}}{\sqrt{\mu_\sigma}})\Big\}\fep\varphi_{a} dvdx\nonumber\\
&\quad+\v^2 \intr\intr \bar{A}(t,x_1,v)\fep\varphi_{a}^2 dvdx+\v^{2}\intr\intr \Gamma(\fep,\fep\varphi_{a}) \fep\varphi_{a} dvdx.
\end{align}

We shall estimate the right hand side of \eqref{4.3} term by term. Firstly we notice that  $\{\pt+v\cdot\nax\}\sqrt{\mu_\sigma}/\sqrt{\mu_\sigma}$ is a cubic polynomial in $v$, then for any $\kappa>0$ and $\delta=\frac1{2(3-\gamma)}$, one has that
\begin{align}\label{4.7}
&\left| \intr\intr \frac{\{\pt+v\cdot\nax\}\sqrt{\mu_\sigma}}{\sqrt{\mu_\sigma}} |\fep\varphi_{a}|^2 dv dx\right|\nonumber\\
&\leq C\intr\intr |\partial_{x_1}(\rho_0,u_0,\theta_0)(t,x_1)|\cdot (1+|v|^2)^{\frac32}\cdot |\fep\varphi_{a}|^2 dv dx\nonumber\\
&=\int\int_{|v|\geq \frac{\kappa}{\ep^{\delta}}}+\int\int_{|v|\leq \frac{\kappa}{\ep^{\delta}}}\nonumber\\
&\leq C \left\{\iint |\partial_{x_1}(\rho_0,u_0,\theta_0)(t,x_1)|^2\cdot (1+|v|^{2})^{3} |\fep\varphi_{a}|^2 I_{\{|v|\geq\kappa/\ep^{\delta}\}}dvdx\right\}^{\frac12}\cdot\|\fep\varphi_{a}\|_{L^2}  \nonumber\\
&\quad +C\|\partial_{x_1}(\rho_0,u_0,\theta_0)\|_{L^{\infty}_{x_{1}}}\cdot\|(1+|v|^{2})^{3/4}\fep\varphi_{a}I_{\{|v|\leq\kappa/\ep^{\delta}\}}\|_{L^2}^{2}\nonumber\\
&\leq Ca^{-2}\|\partial_{x_1}(\rho_0,u_0,\theta_0)\|_{L^{2}_{x_{1}}}\cdot \|h^\v\|_{L^\infty}\cdot \left\{\int_{|v|\geq\frac{\kappa}{\v^{\delta}}} (1+|v|^{2})^{-2\beta+3}  dv\right\}^{\frac12}\|\fep\varphi_{a}\|_{L^2} \nonumber\\
&\quad +C\|\partial_{x_1}(\rho_0,u_0,\theta_0)\|_{L^{\infty}_{x_{1}}}\cdot\|(1+|v|^{2})^{3/4}\fep\varphi_{a}I_{\{|v|\leq\kappa/\ep^{\delta}\}}\|_{L^2}^{2}\nonumber\\
&\leq C_{\kappa} \frac{a^{-2}\v^{2}}{\sqrt{\sigma+t}} \|\hep\|_{L^\infty}\cdot\|\fep\varphi_{a}\|_{L^2}+\frac{C}{\s+t}\|(1+|v|^{2})^{3/4}\mathbf{P}(\fep\varphi_{a})I_{\{|v|\leq\kappa/\ep^{\delta}\}}\|_{L^2}^{2}\nonumber\\
&\quad +\frac{C}{\s+t}\|(1+|v|^{2})^{3/4}\{\mathbf{I-P}\}(\fep\varphi_{a})I_{\{|v|\leq\kappa/\ep^{\delta}\}}\|_{L^2}^{2}\nonumber\\
&\leq C_{\kappa} \frac{a^{-2}\v^{2}}{\sqrt{\sigma+t}} \|\hep\|_{L^\infty}\cdot\|\fep\varphi_{a}\|_{L^2}+\frac{C}{\s+t}\|\fep\varphi_{a}\|_{L^2}^{2}+\frac{C\kappa^{3-\gamma}}{\ep^{\f12}\s}\|\{\mathbf{I-P}\}(\fep\varphi_{a})\|^{2}_{\nu},
	\end{align}
where we have used the fact that $\mu_{M}\leq C\mu_{\sigma}$ (see \eqref{relation of mu and muM}) and
\begin{align}\label{4.8}
|(1+|v|^{2})^{3/2}\fep|&= |(1+|v|^{2})^{-\beta+3/2}\frac{\sqrt{\mu_M}}{\sqrt{\mu_{\sigma}}}\hep|\nonumber\\
&\leq  Ce^{-c_1|v|^2}|(1+|v|^{2})^{-\frac34-2(3-\gamma)}\hep|,
\end{align}
for $\beta\geq \frac{9}{4}+2(3-\gamma)$, where $c_1>0$ is some positive constant depending only on $\theta_M, \theta_-$ and $\theta_+$. \vspace{1.5mm}
	
The appearance of  second term on the right hand side of \eqref{4.3} is mainly due to the cut-off function $\varphi_a$, and it  has not appeared in previous works \cite{Guo Jang Jiang,Guo Jang}. Noting
\begin{align}\nonumber
|v\cdot\nabla_x\varphi_{a}|=a^{-1}\varphi_{a} \frac{|2v\cdot \frac{x}{a}|}{(1-|\frac{x}{a}|^2)^2},
\end{align}
which, together with \eqref{4.8}, yields that
\begin{align}\label{4.9}
&\left|\intr\intr v\cdot\nax\varphi_{a}|\fep|^2\varphi_{a} dx dv\right|\nonumber\\
&\leq \frac{C_\lambda}{a^{1+3\lambda}}\intr\intr|v|\cdot|\fep|^2\varphi_{a}^{2-\lambda}dx dv\nonumber\\
&\leq \frac{C_\lambda}{a^{1+3\lambda}}\|\hep\|_{L^\infty}^{\lambda}\intr\intr |v|\exp{(-c_1 \lambda |v|^2)} |\fep\varphi_{a}|^{2-\lambda}dx dv\nonumber\\
&\leq \frac{C_\lambda}{a^{1+3\lambda}}\|\hep\|_{L^\infty}^{\lambda} \left\{\int_{|x|\leq a}\intr  |v|^{\frac2\lambda}\exp{(-2c_1 |v|^2)}dx dv\right\}^{\frac{\lambda}{2}}\cdot \|\fep\varphi_{a}\|_{L^2}^{2-\lambda}\nonumber\\
&\leq \frac{C_{\lambda}}{a^{1+\frac{3}{2}\lambda}}\|\hep\|_{L^\infty}^{\lambda}\cdot\|\fep\varphi_{a}\|_{L^2}^{2-\lambda}.
\end{align}
where $\lambda\in (0,1)$ is a small constant chosen later.

For the third term on RHS of \eqref{4.3}, we notice that the upper bound of  $F_i$ involving  $\frac{1}{\sigma}$, then for the case $i=1$ we do not have any decay for $\v$. Fortunately, we find that $\Gamma\left(\frac{F_{i}}{\sqrt{\mu_\sigma}},\fep\varphi_{a}\right)+\Gamma\left(\fep\varphi_{a},\frac{F_{i}}{\sqrt{\mu_\sigma}}\right)$ is indeed microscopic part, then one has that
\begin{align}\label{4.10}
&\sum_{i=1}^{5}\ep^{i-1}\intr\intr\left\{\Gamma\left(\frac{F_{i}}{\sqrt{\mu_\sigma}},\fep\varphi_{a}\right)+\Gamma\left(\fep\varphi_{a},\frac{F_{i}}{\sqrt{\mu_\sigma}}\right)\right\}\cdot\fep\varphi_{a}dvdx\nonumber\\
&=\sum_{i=1}^{5}\ep^{i-1}\intr\intr\left\{\Gamma\left(\frac{F_{i}}{\sqrt{\mu_\sigma}},\fep\varphi_{a}\right)+\Gamma\left(\fep\varphi_{a},\frac{F_{i}}{\sqrt{\mu_\sigma}}\right)\right\} \cdot\{\mathbf{I-P}\}\fep\varphi_{a} dvdx\nonumber\\
&\leq \sum_{i=1}^5 (1+t)^{C_0i}\ep^{i-1}\sigma^{-C_0 i}\Big(\|\fep\varphi_a\|_{\n}+\|\fep\varphi_a\|_{L^2}\Big)\cdot\|\{\mathbf{I-P}\}\fep\varphi_a\|_{\n}\nonumber\\
&\leq \sum_{i=1}^5 (1+t)^{C_0i}\ep^{i-1}\s^{-C_0 i}\Big(\|\{\mathbf{I-P}\}(\fep\varphi_a)\|_{\nu}+\|\fep\varphi_a\|_{L^2}\Big)\cdot\|\{\mathbf{I-P}\}\fep\varphi_a\|_{\nu}\nonumber\\
&\leq \Big\{\frac{c_0}{4} +\sum_{i=1}^5 (1+t)^{C_0i}\v^{i}\s^{-C_0i}\Big\}\cdot \frac1{\v}\|\{\mathbf{I-P}\}\fep\varphi_a\|^2_{\nu} \nonumber\\
&\quad+C \sum_{i=1}^5(1+t)^{2C_0i} \ep^{2i-1}\s^{-2C_0 i} \|\fep\varphi_a\|^2_{L^2},
\end{align}
where we have used \eqref{3.26-1} for $F_i, i=1,\cdots,5$ in Theorem \ref{thm1.2}.

For the forth and fifth terms on RHS of \eqref{4.3}, by using \eqref{3.26-1} and \eqref{3.27-1}, it is direct to have that
\begin{align}\label{4.11}
&\left| \v^2 \intr\intr \bar{A}(t,x_1,v)\fep\varphi_{a}^2 dvdx\right|\nonumber\\
&\leq C\v^2\left\{ \intr\intr |\bar{A}(t,x_1,v)\varphi_{a}|^2 dvdx\right\}^{\frac12}\cdot \|\fep\varphi_{a}\|_{L^2}\nonumber\\
&\leq C a^{-\frac32} \v^2 \|\fep\varphi_{a}\|_{L^2} \Big\{\sum_{\substack{i+j\geq 6\\1\leq i,j\leq 5}} (1+t)^{C_0 (i+j)}\cdot \v^{i+j-6}\sigma^{-C_0(i+j)}+(1+t)^{6C_0} \cdot \sigma^{-6C_0}\Big\}\nonumber\\
&\leq C(1+t)^{10C_0}\v^2 \sigma^{-10C_0} \|\fep\varphi_{a}\|_{L^2}.
\end{align}
and
\begin{equation}\label{4.12}
\v^{2} |\la\Gamma(\fep,\fep\varphi_{a}),\fep\varphi_{a}\ra|\leq C\v^{2}\|\hep\|_{L^\infty}\cdot\|\fep\varphi_{a}\|_{L^2}^{2}.
\end{equation}	

Now substituting \eqref{4.7}, \eqref{4.9}-\eqref{4.12} into \eqref{4.3}, one has that
\begin{align}
&\frac{d}{dt}\|\fep\varphi_{a}\|_{L^2}^{2}+\Big\{\frac32c_0-C\kappa^{3-\gamma} -C\v (1+t)^{5C_0}\s^{-5C_0}\Big\}\frac{1}{\ep}\|\{\mathbf{I-P}\}(\fep\varphi_{a})\|_{\nu}^{2}\nonumber\\
&\leq \Big\{C_\kappa a^2\ep^{2}\sigma^{-\f12}\|\hep\|_{L^\infty}+\frac{C_\lambda}{a^{1+\frac{3}{2}\lambda}}\|\hep\|_{L^\infty}^\lambda\nonumber\\
&\qquad\qquad\qquad+C(1+t)^{10C_0} \v \sigma^{-10C_0}+\frac{C}{\sigma+t}\Big\}
\cdot(\|\fep\varphi_{a}\|_{L^2}^{2}+1).\nonumber
\end{align}
Taking $0<C\kappa^{3-\gamma}\leq \frac{c_0}{4}$, and noting $T\leq \v^{-\frac{1}{10C_0}}$ and $\sigma=\v^{\eta}$ with $0<\eta\leq \frac{1}{11C_0}$, one proves \eqref{f L2} by taking $\v\in(0,\v_0)$ with $\v_0$ suitably small. Therefore the proof of Lemma \ref{lem5.1} is completed. $\hfill\Box$

\

\subsection{Weighted $L^\infty$-estimate}
As in \cite{Guo Jang Jiang-1,Guo Jang Jiang}, we denote
\begin{equation}
L_{M}g=-\frac{1}{\smum}\Big\{Q(\mu_\sigma,\smum g)+Q(\smum g,\mu_\sigma)\Big\}=\nu(\mu_{\sigma})g+Kg,\nonumber
\end{equation}
where the frequency $\nu(\mu_{\sigma})$ has been defined in \eqref{2.3} and  $Kg=K_{1}g-K_{2}g$ with
\begin{align}
K_{1}g=&\intr\ints B(\theta)|u-v|^{\gamma}\sqrt{\mu_{M}(u)}\frac{\mu_{\s}(v)}{\sqrt{\mu_{M}(v)}}g(u)dud\omega,
\nonumber\\
K_{2}g=&\intr\ints B(\theta)|u-v|^{\gamma}\mu_\s(u')\frac{\sqrt{\mu_{M}(v')}}{\sqrt{\mu_{M}(v)}}g(v')dud\omega
\nonumber\\
&+\intr\ints B(\theta)|u-v|^{\gamma}\mu_\s(v')\frac{\sqrt{\mu_{M}(u')}}{\sqrt{\mu_{M}(v)}}g(v')dud\omega.\nonumber
\end{align}
Let $0\leq \chi_{m}\leq 1$ be a smooth cut off function, such that for any $m>0$,
\begin{equation}\nonumber
\chi_{m}(s)\equiv 1\ \mbox{for} \  s\leq m, \quad \chi_{m}(s)\equiv0,\  \mbox{for} \  s\geq 2m.
\end{equation}
Then one can define
\begin{align}
K^{m}g=&\intr\ints B(\theta)|u-v|^{\gamma}\chi_{m}(|u-v|)\sqrt{\mu_{M}(u)}\frac{\mu_\s(v)}{\sqrt{\mu_{M}(v)}}g(u)dud\omega\nonumber\\
&+\intr\ints B(\theta)|u-v|^{\gamma}\chi_{m}(|u-v|)\mu_\s(u')\frac{\sqrt{\mu_{M}(v')}}{\sqrt{\mu_{M}(v)}}g(v')dud\omega\nonumber\\
&+\intr\ints B(\theta)|u-v|^{\gamma}\chi_{m}(|u-v|)\mu_\s(v')\frac{\sqrt{\mu_{M}(u')}}{\sqrt{\mu_{M}(v)}}g(v')dud\omega,\nonumber
\end{align}
and  
\begin{equation}\nonumber
K^{c}=K-K^{m}.
\end{equation}

\begin{lemma}[\cite{Guo Jang Jiang,DHWY}]\label{lem5.2-1}
There exists some positive constant  $c>0$, such that
\begin{align}\label{4.16}
|K^mg(v)|\leq Cm^{3+\gamma} e^{-\frac{c}{10}|v|^2} \|g\|_{L^\infty},
\end{align}
and $\displaystyle K^cg(v)=\intr l(v,v')g(v') dv'$ where the kernel $l(v,v')$ satisfies
\begin{equation}\label{4.17}
|l(v,v')|\leq C_{m}\frac{\exp{\{-c|v-v'|^{2}\}}}{|v-v'|(1+|v|+|v'|)^{1-\gamma}}+C|v-v'|^{\gamma} e^{-c|v|^2-c|v'|^2},
\end{equation}
and
\begin{equation}\label{4.17-1}
|l(v,v')|\leq C|v-v'|^{-\frac{3-\gamma}{2}} e^{-c|v-v'|^{2}} e^{-\frac{c||v|^2-|v'|^2|^2}{|v-v'|^2}}+C|v-v'|^{\gamma} e^{-c|v|^2-c|v'|^2}.
\end{equation}
It is worth to point out that the constant $C_m$ is independent of $\sigma$.
\end{lemma}

\begin{lemma}\label{lem5.2}
Let $\eta\leq \frac{1}{40C_0}$, $T= \v^{-\eta}$, $a=\v^{-2\eta}$ and $\sigma=\v^{-\eta}$, then there exists $\ep_0>0$ such that for all $\ep<\ep_0$, $t\in[0,T]$, it holds that
\begin{align}\label{h infty}
&\sup_{0\leq s\leq t}\|\frac{\ep^{3/2}}{a^{3}}\hep(s)\|_{L^\infty}
\leq C\Big\{\|\frac{\ep^{3/2}}{a^{3}}\hep(0)\|_{L^\infty}+C\frac{\ep^{9/2}}{a^{3}}(1+t)^{10C_0}\cdot \sigma^{-10C_0}\Big\}\nonumber\\
&\qquad\qquad\qquad+C\ep^{3/2}a^{3}\sup_{0\leq s\leq t}\|\frac{\ep^{3/2}}{a^{3}}\hep(s)\|^{2}_{L^\infty}+C\sup_{x\in\mathbb{R}^3}\sup_{0\leq s\leq t}\|\fep(s)\varphi_{a}(\cdot-x)\|_{L^2}.
\end{align}
\end{lemma}

\noindent{\bf Proof.}	Letting $K_{w}g\equiv wK(\frac{g}{w})$, it follows from \eqref{Remainder of FEP} and \eqref{def of h} that
\begin{align}
& \pt\hep+v\cdot\nax\hep+\frac{\nu(\mu_\sigma)}{\ep}\hep+\frac{1}{\ep}K_{w}\hep\nonumber\\
&=\ep^{2}\frac{w}{\smum}Q\left(\frac{\hep\smum}{w},\frac{\hep\smum}{w}\right)\nonumber\\
&\quad+\sum_{i=1}^{5}\ep^{i-1}\frac{w}{\smum}\left\{Q\left(F_{i},\frac{\hep\smum}{w}\right) +Q\left(\frac{\hep\smum}{w},F_{i}\right)\right\}+\v^2\tilde{A}(t,x_1,v),\nonumber
\end{align}
with
\begin{align}\nonumber
\tilde{A}(t,x_1,v):=\sum_{\substack{i+j\geq 6\\1\leq i,j\leq 5}}\ep^{i+j-6}\frac{1}{\sqrt{\mu_{M}}}Q(F_{i},F_{j})-\frac{\{\partial_t+v_1 \partial_{x_1}\} F_5}{\sqrt{\mu_{M}}}.
\end{align}
For any $(t,x,v)$, integrating along the backward trajectory, one has that
\begin{align}\label{mild1}
&\hep(t,x,v)\nonumber\\
&=\exp{\left\{-\frac{1}{\ep}\int_{0}^{t}\nu(\tau)d\tau\right\}}\hep(0,x-vt,v)\nonumber\\
&\quad -\frac{1}{\ep}\int_{0}^{t}\exp{\left\{-\frac{1}{\ep}\int_{s}^{t}\nu(\tau)d\tau\right\}} (K_{w}^{m}\hep) (s,x-v(t-s),v)ds\nonumber\\
&\quad -\frac{1}{\ep}\int_{0}^{t}\exp{\left\{-\frac{1}{\ep}\int_{s}^{t}\nu(\tau)d\tau\right\}} (K_{w}^{c}\hep)(s,x-v(t-s),v)ds\nonumber\\
 &\quad +\ep^{2}\int_{0}^{t}\exp{\left\{-\frac{1}{\ep}\int_{s}^{t}\nu(\tau)d\tau\right\}}\left(\frac{w}{\smum}Q\left(\frac{\hep\smum}{w},\frac{\hep\smum}{w}\right)\right)(s,x-v(t-s),v)ds\nonumber\\
&\quad +\int_{0}^{t}\exp{\left\{-\frac{1}{\ep}\int_{s}^{t}\nu(\tau)d\tau\right\}}\left(\sum_{i=1}^{5}\ep^{i-1}\frac{w}{\smum}Q\left(F_{i},\frac{\hep\smum}{w}\right)\right)(s,x-v(t-s),v)ds\nonumber\\
&\quad+\int_{0}^{t}\exp{\left\{-\frac{1}{\ep}\int_{s}^{t}\nu(\tau)d\tau\right\}}\left(\sum_{i=1}^{5}\ep^{i-1}\frac{w}{\smum}Q\left(\frac{\hep\smum}{w},F_{i}\right)\right)(s,x-v(t-s),v)ds\nonumber\\
&\quad +\v^2\int_{0}^{t}\exp{\left\{-\frac{1}{\ep}\int_{s}^{t}\nu(\tau)d\tau\right\}}\tilde{A}(s,x_1-v_1(t-s),v)ds.
\end{align}
It is easy to know that
\begin{align}\label{4.22}
\left|\exp{\left\{-\frac{1}{\ep}\int_{0}^{t}\nu(\tau)d\tau\right\}}\hep(0,x-vt,v)\right|\leq C\|h^\v(0)\|_{L^\infty}.
\end{align}
A direct calculation shows that
\begin{equation}\label{4.23}
\nu(\mu_\s)\sim\nu_{M}(v):=\intr\ints B(v-u,\t)\mu_M(u)d\o du,
\end{equation}
and
\begin{align}\label{4.24}
\int_{0}^{t}\exp{\left\{-\frac{1}{\ep}\int_{s}^{t}\nu(\tau)d\tau\right\}}\nu(\mu_\s)ds &\leq c\int_{0}^{t}\exp{\left\{-\frac{c\nu_{M}(t-s)}{\ep}\right\}}\nu_{M}ds\nonumber\\
&=O(\ep),
\end{align}
where all the constants above are independent of $\sigma$. For the second term on RHS of  \eqref{mild1}, by using \eqref{4.16}, \eqref{4.23} and \eqref{4.24}, it is bounded by
\begin{align}\label{4.25} &\frac{Cm^{3+\gamma}}{\ep}\int_{0}^{t}\exp{\left\{-\frac{1}{\ep}\int_{s}^{t}\nu(\tau)d\tau\right\}}\nu ds\cdot \sup_{0\leq s\leq t}\|\hep(s)\|_{L^\infty} \nonumber\\
&\leq Cm^{3+\gamma}\sup_{0\leq s\leq t}\|\hep(s)\|_{L^\infty}.
\end{align}

\

Since $\mu_{M}\leq C\mu_{\sigma}$, it is easy to know that
\begin{equation}\nonumber
\left|\frac{w}{\smum}Q\left(\frac{\hep\smum}{w},\frac{\hep\smum}{w}\right)\right|\leq C\nu_{M} \|\hep\|_{L^\infty}^{2}\leq  C\nu(\mu_\sigma)\|\hep\|_{L^\infty}^{2},
\end{equation}
the fourth term on RHS of \eqref{mild1} is bounded by
\begin{equation}\label{nonl}
C\ep^{2}\int_{0}^{t}\exp{\left\{-\frac{1}{\ep}\int_{s}^{t}\nu(\tau)d\tau\right\}}\nu(\mu_\sigma)\|\hep(s)\|_{L^\infty}^{2}ds\leq C\ep^{3}\sup_{0\leq s\leq t}\|\hep(s)\|_{L^\infty}^{2}.
\end{equation}

For the fifth and sixth term on RHS of  \eqref{mild1}, it follows from \eqref{3.26-1}  and \eqref{relation of mu and muM} that
\begin{align}
&\quad \left|\sum_{i=1}^{5}\ep^{i-1}\frac{w}{\smum}\left\{Q\left(F_{i},\frac{\hep\smum}{w}\right)+Q\left(\frac{\hep\smum}{w},F_{i}\right)\right\}(t,x,v)\right|\nonumber\\
&\leq C\nu_{M}(v)\|\hep\|_{L^\infty}\left\|\frac{w}{\sqrt{\mu_M}}\sum_{i=1}^{5}\ep^{i-1}F_{i}\right\|_{L^\infty}\nonumber\\
&\leq C\nu_{M}(v)\|\hep\|_{L^\infty} (1+t)^{5C_0}\s^{-5C_0},\nonumber
\end{align}
which yields that the fifth and sixth term on RHS of  \eqref{mild1} are bounded by
\begin{align}\label{linear}
&C\int_{0}^{t}\exp{\left\{-\frac{1}{\ep}\int_{s}^{t}\nu(\tau)d\tau\right\}}\nu_{M}(v)\|\hep(s)\|_{L^\infty}ds\nonumber\\
&\leq C(1+t)^{5C_0} \v\cdot\sigma^{-5C_0}\sup_{0\leq s\leq t}\|\hep(s)\|_{L^\infty}.
\end{align}

For the last term on RHS of  \eqref{mild1},  it follows from \eqref{3.26-1}, \eqref{3.27-1}   and  \eqref{relation of mu and muM},  that
\begin{align}
|\tilde{A}(t,x_1,v)|&\leq C \mu_M(v)^{\alpha-\frac12} \Big\{(1+t)^{10C_0} \sigma^{-10C_0}+C (1+t)^{5C_0} \sigma^{-5C_0-1}\Big\}\nonumber\\
&\leq C\mu_M(v)^{\alpha-\frac12} (1+t)^{10C_0} \sigma^{-10C_0},\nonumber
\end{align}
which, together with \eqref{4.24}, yields  that the last term  on RHS of  \eqref{mild1} is bounded by
\begin{align}\label{4.28}
C\v^3  (1+t)^{10C_0} \sigma^{-10C_0}.
\end{align}

From the definition of $K_{w}^c$ in Lemma \ref{lem5.2-1}, we can bound  the third term on RHS of \eqref{mild1} by
\begin{equation}\label{4.29}
\frac{1}{\ep}\int_{0}^{t}\exp{\left\{-\frac{1}{\ep}\int_{s}^{t}\nu(\tau)d\tau\right\}}\intr|l(v,v')\hep(s,x-v(t-s),v')|dv'ds.
\end{equation}
Using \eqref{mild1} again to \eqref{4.29}, then \eqref{4.29} is bounded by
\begin{align}\label{mild2}
&\frac{1}{\ep}\int_{0}^{t}\exp{\left\{-\frac{1}{\ep}\int_{s}^{t}\nu(\tau)d\tau-\frac{1}{\ep}\int_{0}^{s}\nu(v')(\tau)d\tau\right\}}\intr|l(v,v')|dv'\nonumber\\
&\quad \times |\hep(0,\tilde{x}-v's,v')|ds\nonumber\\
&\quad+\frac{1}{\ep^{2}}\int_{0}^{t}\exp{\left\{-\frac{1}{\ep}\int_{s}^{t}\nu(\tau)d\tau\right\}}\intr|l(v,v')|\int_{0}^{s}\exp{\left\{-\frac{1}{\ep}\int_{s_{1}}^{s}\nu(v')(\tau)d\tau\right\}}\nonumber\\
&\qquad\times|\{K^{m}\hep\}(s_{1},\tilde{x}-v'(s-s_{1}),v')|dv'ds_{1}ds\nonumber\\
&\quad+\frac{1}{\ep^{2}}\int_{0}^{t}\exp{\left\{-\frac{1}{\ep}\int_{s}^{t}\nu(\tau)d\tau\right\}}\intr\intr|l(v,v')l(v',v'')|\nonumber\\
& \qquad\times \int_{0}^{s}\exp{\left\{-\frac{1}{\ep}\int_{s_{1}}^{s}\nu(v')(\tau)d\tau\right\}}|\hep(s_{1},\tilde{x}-v'(s-s_{1}),v'')|dv'dv''ds_{1}ds\nonumber\\
&\quad+\frac{C}{\ep}\int_{0}^{t}\exp{\left\{-\frac{1}{\ep}\int_{s}^{t}\nu(\tau)d\tau\right\}}ds\cdot \intr|l(v,v')|dv'\cdot \{\ep^{3}\sup_{0\leq s\leq t}\|\hep(s)\|_{\infty}^{2}\}\nonumber\\
&\quad+\frac{C}{\ep}\int_{0}^{t}\exp{\left\{-\frac{1}{\ep}\int_{s}^{t}\nu(\tau)d\tau\right\}}ds\cdot \intr|l(v,v')|dv'\nonumber\\
&\qquad\qquad\qquad\qquad\times \{(1+t)^{5C_0} \v\cdot\sigma^{-5C_0}\sup_{0\leq s\leq t}\|\hep(s)\|_{L^\infty}\}\nonumber\\
&\quad+\frac{C}{\ep}\int_{0}^{t}\exp{\left\{-\frac{1}{\ep}\int_{s}^{t}\nu(\tau)d\tau\right\}}ds\cdot \intr|l(v,v')|dv'\cdot \v^3  (1+t)^{10C_0} \sigma^{-10C_0}.
\end{align}
where we have used \eqref{nonl}, \eqref{linear}, \eqref{4.28}, and  denoted $\tilde{x}=x-v(t-s)$ for simplicity of presentation.

It follows from \eqref{4.17} and \eqref{4.17-1} that
\begin{align}\label{4.31-1}
\intr|l(v,v')|dv'\leq
\begin{cases}
C_m (1+|v|^2)^{\frac{\gamma}{2}},\\[2mm]
C(1+|v|)^{-1},
\end{cases}
\end{align}
which yields that  the last three terms and the first term in \eqref{mild2}  are bounded by
\begin{align}\label{4.31}
&C_m\Big\{\|\hep(0)\|_{L^\infty}+\ep^{3}\sup_{0\leq s\leq t}\|\hep(s)\|_{L^\infty}^{2}\nonumber\\
&\qquad+(1+t)^{5C_0} \v\cdot\sigma^{-5C_0}\sup_{0\leq s\leq t}\|\hep(s)\|_{L^\infty}+\v^3  (1+t)^{10C_0} \sigma^{-10C_0}\Big\}.
\end{align}
For the second term in \eqref{mild2}, by using \eqref{4.16}, one can bound it by
\begin{align}\label{4.34}
& \frac{Cm^{3+\gamma}}{\ep^{2}}\sup_{0\leq s\leq t}\|\hep(\tau)\|_{L^\infty}\int_{0}^{t}\exp{\left\{-\frac{\nu_{M}(v)(t-s)}{C\ep}\right\}}\nonumber\\
&\quad \times \intr|l(v,v')|\int_{0}^{s}\exp{\left\{-\frac{\nu_{M}(v')(s-s_{1})}{C\ep}\right\}} e^{-\frac{c}{10}|v'|^2}dv'ds_{1}ds\nonumber\\
&\quad\quad \leq \frac{Cm^{3+\gamma}}{\ep}\sup_{0\leq d\leq t}\|\hep(\tau)\|_{L^\infty}\int_{0}^{t}\exp{\left\{-\frac{\nu_{M}(v)(t-s)}{C\ep}\right\}}\intr|l(v,v')| e^{-\frac{c}{20}|v'|^2} dv'ds\nonumber\\
& \quad\quad \leq Cm^{3+\gamma}\sup_{0\leq \tau\leq t}\|\hep(\tau)\|_{L^\infty}.
\end{align}
	
We now concentrate on the third term in \eqref{mild2}.  As in \cite{Guo Jang Jiang}, we divide it into the following several cases.
	
\noindent{\it Case 1.}  For $|v|\geq N$, by using  $\eqref{4.31-1}_1$, one deduces the following bound:
\begin{align}\label{4.35}
& \frac{C}{\ep^2}\sup_{0\leq s\leq t}\|\hep(s)\|_{L^\infty}\int_{0}^{t}\exp{\left\{-\frac{\nu_{M}(v)(t-s)}{C\ep}\right\}}\intr|l(v,v')|\nonumber\\
&\quad \times \int_{0}^{s}\exp{\left\{-\frac{\nu_{M}(v')(s-s_{1})}{C\ep}\right\}}\intr|l(v',v'')|dv''ds_{1}dv'ds\nonumber\\
&\leq \frac{C_m}{N}\sup_{0\leq s\leq t}\|\hep(s)\|_{L^\infty}.
	\end{align}
	
\noindent{\it Case 2.} For either $|v|\leq N, |v'|\geq 2N$ or $|v'|\leq 2N, |v''|\geq 3N$,  notice that we get either $|v-v'|\geq N$ or $|v'-v''|\geq N$, then either one of the following is valid for some small positive constant $0<c_1\leq \frac{c}{32}$ (where $c>0$ is the one in Lemma \ref{lem5.2-1}):
\begin{equation}\nonumber
\begin{split}
\dis |l(v,v')|&\leq e^{-c_1N^{2}}|l(v,v')e^{c_1|v-v'|^{2}}|,\\[1mm]
\dis |l(v',v'')|&\leq e^{-c_1N^{2}}|l(v',v'')e^{c_1|v'-v''|^{2}}|,
\end{split}
\end{equation}
which, together with \eqref{4.17}, yields that
\begin{equation}\label{4.37}
\begin{split}
\dis \int|l(v,v')e^{c_1|v-v'|^{2}}|dv'\leq C\nu(v),\\
\dis \int|l(v',v'')e^{c_1|v-v'|^{2}}|dv''\leq C\nu(v').
\end{split}
\end{equation}
Hence, for  the case of $|v-v'|\geq N$ or $|v'-v''|\geq N$, it follows from \eqref{4.37} that
\begin{align}\label{expsmall}
& \int_{0}^{t}\int_{0}^{s}\left\{\int\int_{|v|\leq N,|v'|\geq 2N}+\int\int_{|v'|\leq 2N,|v''|\geq 3N}\right\}(\cdots) dv''dv' ds_1ds\nonumber\\
& \quad \leq \frac{C_{m}}{\ep^{2}}e^{-c_1N^2}\sup_{0\leq s\leq t}\|\hep(s)\|_{L^\infty}\int_{0}^{t}\int_{0}^{s}\int|l(v,v')|\exp{\left\{-\frac{\nu_{M}(v)(t-s)}{C\ep}\right\}}\nonumber\\
& \quad\quad \times\exp{\left\{-\frac{\nu_{M}(v')(s-s_{1})}{C\ep}\right\}}\nu_{M}(v')dv'ds_{1}ds\nonumber\\
&\quad \leq C_{m}e^{-c_1N^2}\sup_{0\leq s\leq t}\|\hep(s)\|_{L^\infty}.
\end{align}
	
\noindent{\it Case 3a.} $|v|\leq N, |v'|\leq 2N, |v''|\leq 3N$. In this case, we note $\nu_{M}(v)\geq c_{N}$. Further more, we assume that $s-s_{1}\leq \ep\kappa$ for some small $\kappa>0$ determined later. Then the corresponding part of the third term in \eqref{mild2} is bounded by
\begin{align}\label{intsmall}
&\frac{C}{\ep^{2}}\int_{0}^{t}\int_{s-\ep\kappa}^{s}\exp{\left\{-\frac{c_{N}(t-s)}{\ep}\right\}}\exp{\left\{-\frac{c_{N} (s-s_{1})}{\ep}\right\}}\|\hep(s_{1})\|_{L^\infty}ds_{1}ds\nonumber\\
&\leq C_{N}\sup_{0\leq s\leq t}\{\|\hep(s)\|_{L^\infty}\}\cdot \frac{1}{\ep}\int_{0}^{t}\exp{\left\{-\frac{c_{N}(v)(t-s)}{\ep}\right\}}ds\cdot \int_{s-\ep\kappa}^{s}\frac{1}{\ep}ds_{1}\nonumber\\
&\leq \kappa C_{N}\sup_{0\leq s\leq t}\{\|\hep(s)\|_{L^\infty}\}.
	\end{align}
	
\noindent{\it Case 3b.} $|v|\leq N, |v'|\leq 2N, |v''|\leq 3N$ and $s-s_{1}\geq \ep\kappa$. This is the last remaining case. We can bound the third term in \eqref{mild2} by
\begin{align}\label{main1}
&\frac{C}{\ep^{2}}\int_{0}^{t}\int_{D}\int_{0}^{s-\ep\kappa}\exp{\left\{-\frac{\nu_{M}(v)(t-s)}{C\ep}\right\}}\exp{\left\{-\frac{\nu_{M}(v')(s-s_{1})}{C\ep}\right\}}\nonumber\\[1mm]
&\quad \times |l(v,v')l(v',v'')\hep(s_{1},\tilde{x}-(s-s_{1})v',v'')|ds_{1}dv'dv''ds,
\end{align}
where $D=\{|v'|\leq 2N,|v''|\leq 3N\}$ and $\tilde{x}=x-v(t-s)$. From \eqref{4.17}, it is noted that  $l_{w}(v,v')$ has  possible integrable singularity of $\frac{1}{|v-v'|}$. As in \cite{Guo Jang Jiang}, we choose a  smooth function $l_{N}(v,v')$  with compact support such that
\begin{equation}\label{approximation}
\sup_{|p|\leq 3N}\int_{|v'|\leq 3N}|l_{N}(p,v')-l_{w}(p,v')|dv'\leq \frac{1}{N^{1+|\gamma|}}.
\end{equation}
Splitting
\begin{align}\label{4.41}
& l(v,v')l(v',v'')\nonumber\\
& =\{l(v,v')-l_{N}(v,v')\}l(v',v'')\nonumber\\
&\quad+\{l(v',v'')-l_{N}(v',v'')\}l_{N}(v,v')+l_{N}(v,v')l_{N}(v',v''),
\end{align}
then using $\eqref{4.31-1}_1$, \eqref{approximation} and \eqref{4.41}, we can  bound  \eqref{main1} by
\begin{align}\label{transbefore}
&\quad\frac{C}{\ep^{2}}\int_{0}^{t}\int_{D}\int_{0}^{s-\ep\kappa}\exp{\left\{-\frac{\nu_{M}(v)(t-s)}{C\ep}\right\}}\exp{\left\{-\frac{\nu_{M}(v')(s-s_{1})}{C\ep}\right\}}\nonumber\\
&\quad\quad \times |l_{N}(v,v')l_{N}(v',v'')\hep(s_{1},\tilde{x}-(s-s_{1})v',v'')|ds_{1}dv'dv''ds\nonumber\\
&\quad+\frac{C_m}{N}\sup_{0\leq s\leq t}\{\|\hep(s)\|_{L^\infty}\}.
\end{align}
Since $l_{N}(v,v')l_{N}(v',v'')$ is bounded, we first integrate over $v'$ and make a change of variable $y=\tilde{x}-(s-s_{1})v'$ to get
\begin{align}\label{4.44}
& C_{N}\int_{|v'|\leq 2N}|\hep(s_{1},\tilde{x}-(s-s_{1})v',v'')|dv'\nonumber\\
& \leq C_{N}\int_{|v'|\leq 2N}|\fep(s_{1},\tilde{x}-(s-s_{1})v',v'')|dv'\nonumber\\
& \leq C_{N}\left\{\int_{|v'|\leq 2N}|\fep(s_{1},\tilde{x}-(s-s_{1})v',v'')|^{2}dv'\right\}^{1/2}\nonumber\\
& \leq \frac{C_{N}}{\kappa^{3/2}\ep^{3/2}}\left\{\int_{|y-\tilde{x}|\leq (s-s_{1})2N}|\fep(s_{1},y,v'')|^{2}dy\right\}^{1/2},
\end{align}
where we have used $|\frac{dy}{dv'}|\geq \kappa^{3}\ep^{3}$ as $s-s_{1}\geq \kappa\ep$. Using \eqref{4.44}, we can further bound the first term in \eqref{transbefore} by
\begin{align}\label{4.44-1}
&\frac{C_{N,\kappa}}{\ep^{7/2}}\int_{0}^{t}\int_{0}^{s-\kappa\ep}\exp{\left\{-\frac{c_N(t-s)}{\ep}\right\}}\exp{\left\{-\frac{c_{N}(s-s_{1})}{\ep}\right\}}\nonumber\\
& \quad\quad\quad\quad\quad\quad \times\int_{|v''|\leq 3N}\left\{\int_{|y-\tilde{x}|\leq 2N(s-s_1)}|\hep(s_{1},y,v'')|^{2}dy\right\}^{1/2}dv''ds_{1}ds\nonumber\\
& \leq \frac{C_{N,\kappa} a^3}{\ep^{7/2}}\int_{0}^{t}\int_{0}^{s-\kappa\ep}\exp{\left\{-\frac{c_N(t-s)}{\ep}\right\}}\exp{\left\{-\frac{c_N(s-s_{1})}{\ep}\right\}}\nonumber\\
& \quad\quad\quad\quad\quad\quad \times\left\{\int_{|v''|\leq 3N}\int_{|y-x|\leq 2Nt}|\fep(s_{1},y,v'')\varphi_{a}|^{2}dydv''\right\}^{1/2}ds_{1}ds\nonumber\\
& \leq \frac{C_{N,\kappa}a^{3}}{\ep^{3/2}}\sup_{0\leq s\leq t}\|\fep(s)\varphi_{a}(\cdot-x)\|_{L^2},
\end{align}
where we have chosen $a$ to be a positive constant such that  $a\geq 4N(t+1)$.

Collecting all the above terms and multiplying them with $\displaystyle\frac{\ep^{\frac32}}{a^{3}}$,  for any $\kappa>0$ and large $N>0$, then one obtains that
\begin{align}
&\sup_{0\leq s\leq t}\{\|\frac{\ep^{3/2}}{a^{3}}\hep(s)\|_{L^\infty}\}\nonumber\\
&\leq C_m\Big\{\|\frac{\ep^{3/2}}{a^{3}}\hep(0)\|_{L^\infty}+\frac{\ep^{9/2}}{a^{3}}(1+t)^{10C_0}\cdot \sigma^{-10C_0}+\ep^{3/2}a^{3}\|\frac{\ep^{3/2}}{a^{3}}\hep(s)\|^{2}_{L^\infty}\Big\}\nonumber\\
&\quad +C\left\{m^{3+\gamma}+\kappa\cdot C_{N}+C_{m}\big[\frac1N+(1+t)^{5C_0}\v\cdot\sigma^{-5C_0}\big]\right\}\sup_{0\leq s\leq t}\|\frac{\ep^{3/2}}{a^{3}}\hep(s)\|_{L^\infty}\nonumber\\
& \quad+C_{N,\kappa}\sup_{x\in\mathbb{R}^3}\sup_{0\leq s\leq t}\|\fep(s)\varphi_{a}(\cdot-x)\|_{L^2}.\nonumber
\end{align}
Noting $t\in[0,T]$, $T= \v^{\eta}$,$\sigma=\v^{\eta}$, and $ a=\v^{-2\eta}$ with   $\eta\leq \frac{1}{40C_0}$, first choosing $m$ small, then $N$ large enough, and then letting $\kappa$ small, and finally $\v\leq \v_0$ with $\v_0$ small enough so that
$$C\left\{m^{3+\gamma}+\kappa\cdot C_{N}+C_{m}\big[\frac1N+(1+t)^{5C_0}\v\cdot\sigma^{-5C_0}\big]\right\}\leq \frac12,$$
 thus we deduce
\begin{align}\nonumber
&\sup_{0\leq s\leq t}\|\frac{\ep^{3/2}}{a^{3}}\hep(s)\|_{L^\infty}
\leq C\Big\{\|\frac{\ep^{3/2}}{a^{3}}\hep(0)\|_{\infty}+C\frac{\ep^{9/2}}{a^{3}}(1+t)^{10C_0}\cdot \sigma^{-10C_0}\Big\}\nonumber\\
&\qquad\qquad\qquad+C\ep^{3/2}a^{3}\sup_{0\leq s\leq t}\|\frac{\ep^{3/2}}{a^{3}}\hep(s)\|^{2}_{L^\infty}+C\sup_{x\in\mathbb{R}^3}\sup_{0\leq s\leq t}\|\fep(s)\varphi_{a}(\cdot-x)\|_{L^2}.\nonumber
\end{align}
Therefore the proof of 	Lemma \ref{lem5.2} is completed. $\hfill\Box$	
	
	
\subsection{Proof of Theorem \ref{theorem}} Throughout this subsection, we assume that  $T= \v^{-\eta}$, $a=\v^{-2\eta}$ and $\sigma=\v^{\eta}$ where we choose  $\eta:= \min\{\frac1{40C_0},\frac1{100\tilde{C}_1}\}$, and $C_0\geq 1$ and $\tilde{C}_1\geq1$ are the constants determined in Theorem \ref{thm1.2} and Lemma \ref{lem5.1}, respectively.

Now we make the {\it a priori} assumption
\begin{equation}\label{4.46}
\sup_{0\leq t\leq T}\v^{\frac14}\|\frac{\v^{\frac{3}{2}}}{a^3} h^\v(t)\|_{L^\infty}\leq 1,
\end{equation}
then, by taking $\lambda=\frac{1}{21}\eta$, it follows from \eqref{4.46} that
\begin{align}\label{4.47}
&\tilde{C}_1 a^2\ep^{2}\sigma^{-\f12}\|\hep(t)\|_{L^\infty}+\frac{C_\lambda}{a^{1+\frac{3}{2}\lambda}}\|\hep(t)\|_{L^\infty}^\lambda+\tilde{C}_1(1+t)^{10C_0} \v \sigma^{-10C_0}+\frac{\tilde{C}_1}{\sigma+t}\nonumber\\
&\leq 2\tilde{C}_1\v^{2\eta}+C_\eta \cdot \v^{\frac32\eta}+\frac{\tilde{C}_1}{\sigma+t}\nonumber\\
&\leq \frac{4\tilde{C}_1}{\sigma+t},
\end{align}
provided $\v\in(0,\v_0)$ with $\v_0>0$ further small such that  $C_\eta   \v_0^{\frac12\eta}\leq \tilde{C}_1$. Now it follows from \eqref{f L2} and \eqref{4.47} that
\begin{align}\label{4.48}
&\frac{d}{dt}\|\fep(t)\varphi_{a}(\cdot-x_0)\|_{L^2}^{2}\leq  \frac{4\tilde{C}_1}{\sigma+t}
\cdot(\|\fep(t)\varphi_{a}(\cdot-x_0)\|_{L^2}^{2}+1),
\end{align}
which yields immediately that, for $t\in[0,T]$,
\begin{align}\label{4.49}
(\|\fep(t)\varphi_{a}(\cdot-x_0)\|_{L^2}^{2}+1)
&\leq (\|\fep(0)\varphi_{a}(\cdot-x_0)\|_{L^2}^{2}+1)\cdot\left(\frac{\sigma+t}{\sigma} \right)^{4\tilde{C}_1}\nonumber\\
&\leq  (\sup_{x_0\in\mathbb{R}^3}\|\fep(0)\varphi_{a}(\cdot-x_0)\|_{L^2}^{2}+1)\cdot \v^{-8\tilde{C}_1\eta}\nonumber\\
&\leq C\v^{-\frac14-8\tilde{C}_1\eta},
\end{align}
where we have used the initial condition $\sup_{x_0\in\mathbb{R}^3}\|\fep(0)\varphi_{a}(\cdot-x_0)\|_{L^2}^{2}\lesssim \v^{-\frac18}$.
Substituting \eqref{4.49} into \eqref{h infty} and noting \eqref{4.46}, one has that
\begin{align}\label{4.50}
\sup_{0\leq s\leq t}\|\frac{\ep^{3/2}}{a^{3}}\hep(s)\|_{L^\infty}
&\leq C\Big\{\|\frac{\ep^{3/2}}{a^{3}}\hep(0)\|_{L^\infty}+\frac{\ep^{9/2}}{a^{3}}(1+t)^{10C_0}\cdot \sigma^{-10C_0}+\v^{-\frac18-4\tilde{C}_1\eta}\Big\}\nonumber\\
&\leq C\Big\{\|\frac{\ep^{3/2}}{a^{3}}\hep(0)\|_{L^\infty}+\ep^{4}+\v^{-\frac18-4\tilde{C}_1\eta}\Big\}.
\end{align}
Hence noting $\eta=\min\{\frac1{40C_0},\frac1{100\tilde{C}_1}\}$ and using \eqref{4.50}, one obtains that
\begin{align}\label{4.51}
\v^{\frac14}\|\frac{\v^{\frac{3}{2}}}{a^3} h^\v(t)\|_{L^\infty}
&\leq C\v^{\frac1{4}}\Big\{\|\frac{\ep^{3/2}}{a^{3}}\hep(0)\|_{L^\infty}+\v^4\Big\}
+C \v^{\frac18-4\tilde{C}_1\eta}\nonumber\\
&\leq C\v^{\frac1{4}}\Big\{\|\frac{\ep^{3/2}}{a^{3}}\hep(0)\|_{L^\infty}+\v^4\Big\}
+C \v^{\frac{17}{200}}\nonumber\\
&\leq \frac12,
\end{align}
where we have used the initial condition $\|\frac{\ep^{3/2}}{a^{3}}\hep(0)\|_{L^\infty}\lesssim 1$, and $\v_0$ been chosen further small.  In light of \eqref{4.51}, the {\it a priori} assumption \eqref{4.46} will be closed by a continuity argument.

Finally, combining \eqref{4.49} and \eqref{4.51}, we proved \eqref{1.37} and \eqref{1.38}. Therefore the proof of Theorem \ref{theorem} is completed. $\hfill\Box$

\vspace{1.5mm}

	
\section{Appendix}

\noindent{\bf Proof of Lemma \ref{lem4.1}.}
1) 
We need only to show that for any $g_1\in B_p, g_2\in B_q$, $g_1g_2\in B_{p+q}$. Since $B_i$ is linear space,  without loss of generality, we assume that $g_1$ and $g_2$ are the base of $B_p$ and  $B_q$, respectively. It follows from \eqref{def of B} that
\begin{align}\nonumber
g_1g_2=\prod_{i=1}^{p}\pa_\tau^{m_i-1}J_\tau\cdot\prod_{i=1}^{q}\pa_\tau^{\tilde{m}_i-1}J_\tau=\prod_{i=1}^{p+q}\pa_\tau^{\bar{m}_i-1}J_\tau
\end{align}
where $\bar{m}_i=m_i, 1\leq i\leq p$, $\bar{m}_i=\tilde{m}_{i-p}, p+1\leq i\leq p+q$. Hence we have proved $g_1g_2\in B_{p+q}$.\vspace{1.5mm}

2) Let $g$ be any base of $B_n$. Applying the Leibnitz rule, one has that
\begin{align}\nonumber
\pa_\tau g&=\pa_\tau\left(\prod_{i=1}^{n}\pa_\tau^{m_i-1}J_\tau\right)=\sum_{j=1}^n\prod_{i=1}^{n}\pa_\tau^{m_i-1+\delta_{ij}}J_\tau
\end{align}
for each $j$, let $\bar{m}_{i}=m_i+\delta_{ij}, 1\leq i\leq n$ and  $\bar{m}_{n+1}=0$, then $|\bar{m}(n+1)|=\sum_{i=1}^{n}(m_i+\delta_{ij})+0=n+1$, so $\pa_\tau g\in B_{n+1}$. Then it is direct to show that $\pa_\tau (b_nf)=b_{n+1}f+b_n\pa_\tau f$;

3) It follows from \eqref{2.8} that
\begin{align}
|b_k|&\leq C\max_{|m(k)|=k}\left|\prod_{i=1}^k\pa_\tau^{m_i-1}J_\tau\right|
\leq C\max_{|m(k)|=k}\left(\prod_{\substack{i=1\\m_i\neq 0}}^k|\pa_\tau^{m_i-1} J_\tau|\right)\nonumber\\
&\leq C\max_{|m(k)|=k}\left(\prod_{\substack{i=1\\m_i\neq 0}}^k C_{m_i}\frac{\s^{-m_i+1}}{\sigma+t}(1+|v|)^2\right)\leq C\s^{-k}(1+|v|)^{2k}.\nonumber
\end{align}
Therefore the proof of Lemma \ref{lem4.1} is completed.  $\hfill\Box$

\

\noindent{\bf Proof of Lemma \ref{prop of A}.}  Firstly we prove \eqref{prop of A_k}. Noting  $A_1f=\pa_\tau f +\f12 fJ_\tau$, we know that \eqref{prop of A_k} holds for $k=1$. For simplicity, we may denote $A_1f$ with $\pa_\tau f+b_1 f$ instead of $\pa_\tau f +\f12 J_\tau f$. Assume \eqref{prop of A_k} holds for $k=1,...,n$, by the \eqref{def of A}, we have
\begin{align}
A_{n+1}f
&=\pa_\tau\left(\sum_{i=0}^nb_i\pa_\tau^{n-i}f\right)+b_1\left(\sum_{i=0}^nb_i\pa_\tau^{n-i}f\right)\nonumber\\
&=\sum_{i=0}^nb_i\pa_\tau^{n+1-i}f+\sum_{i=0}^k\pa_\tau b_i\pa_\tau^{k-i}f+\sum_{i=0}^kb_{i+1}\pa_\tau^{k-i}f\nonumber\\
&=\sum_{i=0}^nb_i\pa_\tau^{n+1-i}f+\sum_{i=0}^kb_{i+1}\pa_\tau^{(k+1)-(i+1)}f=\sum_{i=0}^{n+1}b_i\pa_\tau^{k+1-i}f,\nonumber
\end{align}
which means \eqref{prop of A_k} holds for $k=n+1$. \vspace{1.5mm}

To prove \eqref{A and G},  a direct calculation shows that
\begin{align}\label{5.2}
A_1\circ \Gamma_{i,j}(f)
&=\pa_\tau\left(\mu_{\s}^{-\f12}[Q(b_i\mu_{\s}^{\f12}\mathbf{L}^{-1}f,b_j\mu_{\s})+Q(b_j\mu_{\s}, b_i\mu_{\s}^{\f12}\mathbf{L}^{-1}f)]\right)+b_1\Gamma_{i,j}(f)\nonumber\\
&= b_1\mu_{\s}^{-\f12}[Q(b_i\mu_{\s}^{\f12}\mathbf{L}^{-1}f,b_j\mu_{\s})+Q(b_j\mu_{\s}, b_i\mu_{\s}^{\f12}\mathbf{L}^{-1}f)]\nonumber\\
&\quad +\mu_{\s}^{-\f12}[Q(\pa_\tau b_i\mu_{\s}^{\f12}\mathbf{L}^{-1}f,b_j\mu_{\s})+Q(b_j\mu_{\s}, \pa_\tau b_i\mu_{\s}^{\f12}\mathbf{L}^{-1}f)]\nonumber\\
&\quad +\mu_{\s}^{-\f12}[Q(b_ib_1\mu_{\s}^{\f12}\mathbf{L}^{-1}f,b_j\mu_{\s})+Q(b_j\mu_{\s}, b_ib_1\mu_{\s}^{\f12} \mathbf{L}^{-1}f)]\nonumber\\
&\quad +\mu_{\s}^{-\f12}[Q(b_i\mu_{\s}^{\f12}\pa_\tau \mathbf{L}^{-1}f,b_j\mu_{\s})+Q(b_j\mu_{\s}, b_i\mu_{\s}^{\f12}\pa_\tau\mathbf{L}^{-1}f)]\nonumber\\
&\quad +\mu_{\s}^{-\f12}[Q(b_i\mu_{\s}^{\f12}\mathbf{L}^{-1}f,\pa_\tau b_j\mu_{\s})+Q(\pa_\tau b_j\mu_{\s}, b_i\mu_{\s}^{\f12}\mathbf{L}^{-1}f)]\nonumber\\
&\quad +\mu_{\s}^{-\f12}[Q(b_i\mu_{\s}^{\f12}\mathbf{L}^{-1}f,b_jb_1\mu_{\s})+Q(b_jb_1\mu_{\s}, b_i\mu_{\s}^{\f12}\mathbf{L}^{-1}f)]\nonumber\\
&=b_1\Gamma_{i,j}(f)+\Gamma_{i+1,j}(f)+\Gamma_{i,j+1}(f)\nonumber\\
&\qquad+\m^{-\f12}[Q(b_i\m^{\f12}\pa_\tau \mathbf{L}^{-1}f,b_j\m)+Q(b_j\m, b_i\m^{\f12}\pa_\tau \mathbf{L}^{-1}f)].
\end{align}
For $f\in \mathcal{N}^{\bot}$, it is noted that
\begin{equation}\nonumber
A_1f=\pa_\tau f+\f12 J_\tau f \in \mathcal{N}^{\bot},
\end{equation}
which, together with a direct calculations, yields that
\begin{align}\label{5.1}
\pa_\tau \mathbf{L}^{-1}f
&=\mathbf{L}^{-1}[A_1f]+b_1\mathbf{L}^{-1}f+\mathbf{L}^{-1}[\Gamma_{0,1}f].
\end{align}
Substituting \eqref{5.1} into \eqref{5.2}, one obtains that
\begin{align}
&A_1\circ \Gamma_{i,j}(f)\nonumber\\
&=b_1\Gamma_{i,j}(f)+\Gamma_{i+1,j}(f)+\Gamma_{i,j+1}(f)\nonumber\\
&\quad +\mu_{\s}^{-\f12}[Q(b_i\mu_{\s}^{\f12}\mathbf{L}^{-1}[A_1f],b_j\mu_{\s})+Q(B_j\mu_{\s}, b_i\mu_{\s}^{\f12}\mathbf{L}^{-1}[A_1f])]\nonumber\\
&\quad +\mu_{\s}^{-\f12}[Q(b_i\mu_{\s}^{\f12}b_1\mathbf{L}^{-1}f,b_j\mu_{\s})+Q(b_j\mu_{\s}, b_i\mu_{\s}^{\f12}b_1\mathbf{L}^{-1}f)]\nonumber\\
&\quad +\mu_{\s}^{-\f12}[Q(b_i\mu_{\s}^{\f12}\mathbf{L}^{-1}[\Gamma_{0,1}f],b_j\mu_{\s})+Q(b_j\mu_{\s}, b_i\mu_{\s}^{\f12}\mathbf{L}^{-1}[\Gamma_{0,1}f])]\nonumber\\
&=b_1\Gamma_{i,j}(f)+\Gamma_{i+1,j}(f)+\Gamma_{i,j+1}(f)+\Gamma_{i,j}\circ A_1(f)+\Gamma_{i,j}\circ \Gamma_{0,1}(f),\nonumber
\end{align}
which proves \eqref{A and G}. Therefore the proof of Lemma  \ref{prop of A} is completed. $\hfill\Box$

\

\noindent{\bf Proof of Lemma \ref{derivatives of \L}.} Recall the notation $N_s(i,j,l)$ in Lemma \ref{derivatives of \L}, and we write $\sum_{(i_m,j_m,l_m)\in N_s(i,j,l)}$ to be $\sum_{N_s(i,j,l)}$ for simplicity of presentation. We shall use induction argument to prove this lemma. Noting \eqref{5.1}, we know that \eqref{n-th derivates of L^-1} holds for $n=0,1$. And we assume that \eqref{n-th derivates of L^-1} holds for  $n$-th derivatives of $\mathbf{L}^{-1}f$. Next, we shall consider the $n+1$-th derivatives of $\mathbf{L}^{-1}f$. By using \eqref{5.1}, a direct calculation shows that
	\begin{align}
&\pa_\tau^{n+1}\mathbf{L}^{-1} f=\pa_\tau(\pa_\tau^n\mathbf{L}^{-1} f)\nonumber\\
&=\sum_{\substack{r+k=n\\r,k\geq0}}\sum_{\substack{s+p=k\\s,p\geq0}}\pa_\tau b_r\mathbf{L}^{-1}[\sum_{\substack{i+j+l=s\\i,j,l\geq0}
}\sum_{N_s(i,j,l)}\left(b_{i_1}\Gamma_{j_1,l_1}\right)\circ\cdot\cdot\cdot\circ\left(b_{i_s}\Gamma_{j_s,l_s}\right)\circ A_pf]\nonumber\\
&\quad +\sum_{\substack{r+k=n\\r,k\geq0}}\sum_{\substack{s+p=k\\s,p\geq0}}b_r\pa_\tau\mathbf{L}^{-1}[\sum_{\substack{i+j+l=s\\i,j,l\geq0}
}\sum_{N_s(i,j,l)}\left(b_{i_1}\Gamma_{j_1,l_1}\right)\circ\cdot\cdot\cdot\circ\left(b_{i_s}\Gamma_{j_s,l_s}\right)\circ A_pf]\nonumber\\
&=\sum_{\substack{r+k=n\\r,k\geq0}}\sum_{\substack{s+p=k\\s,p\geq0}}b_{r+1}\mathbf{L}^{-1}[\sum_{\substack{i+j+l=s\\i,j,l\geq0}
}\sum_{N_s(i,j,l)}\left(b_{i_1}\Gamma_{j_1,l_1}\right)\circ\cdot\cdot\cdot\circ\left(b_{i_s}\Gamma_{j_s,l_s}\right)\circ A_pf]\nonumber\\
&\quad +\sum_{\substack{r+k=n\\r,k\geq0}}\sum_{\substack{s+p=k\\s,p\geq0}}b_r\mathbf{L}^{-1}[\sum_{\substack{i+j+l=s\\i,j,l\geq0}
}\sum_{N_s(i,j,l)}\Gamma_{0,1}\circ\left(b_{i_1}\Gamma_{j_1,l_1}\right)\circ\cdot\cdot\cdot\circ\left(b_{i_s}\Gamma_{j_s,l_s}\right)\circ A_pf]\nonumber\\
&\quad +\sum_{\substack{r+k=n\\r,k\geq0}}\sum_{\substack{s+p=k\\s,p\geq0}}b_r\mathbf{L}^{-1}[\sum_{\substack{i+j+l=s\\i,j,l\geq0}
}\sum_{N_s(i,j,l)}A_1\circ\left(b_{i_1}\Gamma_{j_1,l_1}\right)\circ\cdot\cdot\cdot\circ\left(b_{i_s}\Gamma_{j_s,l_s}\right)\circ A_pf].\label{5.4}
\end{align}
To deal with the last term of \eqref{5.4}, by using \eqref{A and G}, one has that
\begin{align}\label{5.5}
&A_1\circ \left(b_{i_m}\Gamma_{j_m,l_m}\right)(f)\nonumber\\
&=\pa_\tau b_{i_m}\Gamma_{j_m,l_m}(f)+b_{i_m}\pa_\tau\Gamma_{j_m,l_m}(f)+b_{i_m}b_1\Gamma_{j_m,l_m}(f)\nonumber\\
&=b_{i_m+1}\Gamma_{j_m,l_m}(f)+b_{i_m}A_1\circ\Gamma_{j_m,l_m}(f)\nonumber\\
&=\left(b_{i_m}\Gamma_{j_m,l_m}\right)\circ A_1(f)+b_{i_m+1}\Gamma_{j_m,l_m}(f)+b_{i_m}\Gamma_{j_m+1,l_m}(f)\nonumber\\
&\quad+b_{i_m}\Gamma_{j_m,l_m+1}(f)+\left(b_{i_m}\Gamma_{j_m,l_m}\right)
\circ\Gamma_{0,1}(f).
\end{align}
which yields immediately that
\begin{align}\label{5.6}
&\sum_{\substack{r+k=n\\r,k\geq0}}\sum_{\substack{s+p=k\\s,p\geq0}}b_r\mathbf{L}^{-1}[\sum_{\substack{i+j+l=s\\i,j,l\geq0}
}\sum_{N_s(i,j,l)}A_1\circ\left(b_{i_1}\Gamma_{j_1,l_1}\right)\circ\cdot\cdot\cdot\circ\left(b_{i_s}\Gamma_{j_s,l_s}\right)\circ A_pf]\nonumber\\
&=\sum_{\substack{r+k=n\\r,k\geq0}}\sum_{\substack{s+p=k\\s,p\geq0}}b_r\mathbf{L}^{-1}[\sum_{\substack{i+j+l=s\\i,j,l\geq0}
}\sum_{N_s(i,j,l)}\left(b_{i_1+1}\Gamma_{j_1,l_1}\right)\circ\cdot\cdot\cdot\circ\left(b_{i_s}\Gamma_{j_s,l_s}\right)\circ A_pf]\nonumber\\
&\quad +\sum_{\substack{r+k=n\\r,k\geq0}}\sum_{\substack{s+p=k\\s,p\geq0}}b_r\mathbf{L}^{-1}[\sum_{\substack{i+j+l=s\\i,j,l\geq0}
}\sum_{N_s(i,j,l)}\left(b_{i_1}\Gamma_{j_1+1,l_1}\right)\circ\cdot\cdot\cdot\circ\left(b_{i_s}\Gamma_{j_s,l_s}\right)\circ A_pf]\nonumber\\
&\quad +\sum_{\substack{r+k=n\\r,k\geq0}}\sum_{\substack{s+p=k\\s,p\geq0}}b_r\mathbf{L}^{-1}[\sum_{\substack{i+j+l=s\\i,j,l\geq0}
}\sum_{N_s(i,j,l)}\left(b_{i_1}\Gamma_{j_1,l_1+1}\right)\circ\cdot\cdot\cdot\circ\left(b_{i_s}\Gamma_{j_s,l_s}\right)\circ A_pf]\nonumber\\
&\quad +\sum_{\substack{r+k=n\\r,k\geq0}}\sum_{\substack{s+p=k\\s,p\geq0}}b_r\mathbf{L}^{-1}[\sum_{\substack{i+j+l=s\\i,j,l\geq0}
}\sum_{
N_s(i,j,l)}\left(b_{i_1}\Gamma_{j_1,l_1}\right)\circ\Gamma_{0,1}\circ\left(b_{i_2}\Gamma_{j_2,l_2}\right)\circ\cdot\cdot\cdot\circ\left(b_{i_s}\Gamma_{j_s,l_s}\right)\circ A_pf]\nonumber\\
&\quad +\sum_{\substack{r+k=n\\r,k\geq0}}\sum_{\substack{s+p=k\\s,p\geq0}}b_r\mathbf{L}^{-1}[\sum_{\substack{i+j+l=s\\i,j,l\geq0}
}\sum_{N_s(i,j,l)}\left(b_{i_1}\Gamma_{j_1,l_1}\right)\circ A_1\circ\left(b_{i_2}\Gamma_{j_2,l_2}\right)\circ\cdot\cdot\cdot\circ\left(b_{i_s}\Gamma_{j_s,l_s}\right)\circ A_pf].
\end{align}
Again substituting \eqref{5.5} to the last term of \eqref{5.6} until $A_1$ applying to the operator $A_p(f)$, and  changing $(i_m,j_m,l_m)$ to $(i_m+1,j_m,l_m)$, $(i_m,j_m+1,l_m)$ and $(i_m,j_m,l_m+1)$, or add the operator   $\Gamma_{0,1}$ behind
$b_{i_m}\Gamma_{j_m,l_m}$ for each $m=2,...,s$  as in \eqref{5.5}, then noting $A_1\circ A_pf=A_{p+1}f$, one can finally prove that
\begin{align}
&\pa_\tau^{n+1} (\mathbf{L}^{-1} f)\nonumber\\
&=\sum_{\substack{r+k=n+1\\r,k\geq0}}\sum_{\substack{s+p=k\\s,p\geq0}}b_r\mathbf{L}^{-1}[\sum_{\substack{i+j+l=s\\i,j,l\geq0}
}\sum_{N_s(i,j,l)}\left(b_{i_1}\Gamma_{j_1,l_1}\right)\circ\cdot\cdot\cdot\circ\left(b_{i_s}\Gamma_{j_s,l_s}\right)\circ A_pf].\nonumber
\end{align}
Thus we complete the proof of Lemma \ref{derivatives of \L}.  $\hfill\Box$
	
\

\noindent{\bf Proof of Lemma \ref{estimate of Gamma}.}  From the definition of \eqref{2.9}, one has that
\begin{align}
\Gamma_{i,j}(f)
&=\f{1}{\sqrt{\mu_\s(v)}}\int\int_{\mathbb{R}^3\times\mathbb{S}^2}B(v-u,\t)\Big[b_i(v')\sqrt{\mu_\s(v')}\mathbf{L}^{-1} f(v')b_j(u')\m_\s(u')\nonumber\\
&\qquad\qquad\qquad\qquad\qquad\qquad-b_i(v)\sqrt{\m_\s(v)}\mathbf{L}^{-1} f(v)b_j(u)\m_\s(u)\Big]d\o du\nonumber\\
&\quad+\f{1}{\sqrt{\m_\s(v)}}\int\int_{\mathbb{R}^3\times\mathbb{S}^2}B(v-u,\t)\Big[b_i(u')\sqrt{\mu_\s(u')}\mathbf{L}^{-1} f(u')b_j(v')\mu_\s(v')\nonumber\\
&\qquad\qquad\qquad\qquad\qquad\qquad-b_i(u)\sqrt{\mu_\s(u)}\mathbf{L}^{-1} f(u)b_j(v)\mu_\s(v)\Big]d\o du\nonumber\\
&:=I+II.\nonumber
\end{align}

Since $\mathbf{L}^{-1}$ preserves decay in $v$ (see \cite{Caflish}), one has that
\begin{equation}\label{2.15}
|\mathbf{L}^{-1}f|\leq C S(t,x)(1+|v|)^m\sqrt{\mu_{\s}(v)},
\end{equation}
which implies that
\begin{align}
|I|&\leq \f{1}{\sqrt{\mu_\s(v)}}\int\int_{\mathbb{R}^3\times\mathbb{S}^2}|B(v-u,\t)|\cdot|b_i(v')\sqrt{\mu_\s(v')}\mathbf{L}^{-1} f(v')b_j(u')\mu_\s(u')|d\o du\nonumber\\
&\quad + \f{1}{\sqrt{\mu_\s(v)}}\int\int_{\mathbb{R}^3\times\mathbb{S}^2}|B(v-u,\t)|\cdot|b_i(v)\sqrt{\mu_\s(v)}\mathbf{L}^{-1} f(v)b_j(u)\mu_\s(u)|d\o du\nonumber\\
&:=I_1+I_2.\nonumber
\end{align}
Using \eqref{2.12}, one has that
\begin{align}
I_2&\leq\f{C\sigma^{-i-j}}{\sqrt{\mu_\s(v)}}S(t,x)\int\int_{\mathbb{R}^3\times\mathbb{S}^2}|B(v-u,\t)|(1+|v|)^{2i+m}(1+|u|)^{2j}\mu_\s(v)\mu_\s(u)d\o du\nonumber\\
&\leq C_{i,j}S(t,x)\s^{-(i+j)}(1+|v|)^{m+2i+\gamma}\sqrt{\mu_{\s}(v)}.\nonumber
\end{align}
For $I_1$, noting $|v'|\lesssim |v|+|u|,\  |u'|\lesssim |v|+|u|$ and $\m_\s(v')\m_\s(u')=\m_\s(v)\m_\s(u)$, one can obtain
\begin{align}
I_1&\leq \f{C\sigma^{-i-j}}{\sqrt{\m_\s(v)}}S(t,x)\int\int_{\mathbb{R}^3\times\mathbb{S}^2}|B(v-u,\t)|(1+|v'|)^{2i+m}(1+|u'|)^{2j}\m_\s(v')\m_\s(u')d\o du\nonumber\\
&\leq C_{i,j}\sigma^{-(i+j)}S(t,x)\sqrt{\m_\s(v)}\int\int_{\mathbb{R}^3\times\mathbb{S}^2}|B(v-u,\t)|(1+|v|+|u|)^{m+2i+2j}\m_\s(u)d\o du\nonumber\\
&\leq C_{i,j}\sigma^{-(i+j)}S(t,x)(1+|v|)^{m+2i+2j+\gamma}\sqrt{\m_{\s}}_\s(v).\nonumber
\end{align}
Thus combining the above estimates, one gets that
\begin{equation}\label{2.17}
|I|\leq C_{i,j}\s^{-(i+j)}S(t,x)(1+|v|)^{m+2i+2j+\gamma}\sqrt{\mu_\s(v)}
\end{equation}
By similar arguments, one can obtain
\begin{equation}\nonumber
|II|\leq C_{i,j,\w}\s^{-(i+j)}S(t,x)(1+|v|)^{m+2i+2j+\gamma}\sqrt{\mu_\s(v)},
\end{equation}
which, together with \eqref{2.17}, yields \eqref{2.14}. Therefore the proof of Lemma \ref{estimate of Gamma} is completed. $\hfill\Box$

\
	
\noindent{\bf Acknowledgments.} The research of  Yong Wang is partially supported by National Natural Sciences Foundation of China No. 11771429, 11671237 and  11688101.


\begin{thebibliography}{99}
		


\bibitem{Bardos} C. Bardos, F. Golse, C.D.  Levermore,  Fluid dynamic limits of kinetic equations. I. Formal derivations. J. Statist. Phys. 63 (1991), no. 1-2, 323-344.

\bibitem{Bardos-2} C. Bardos, F. Golse, C.D.  Levermore, Fluid dynamic limits of kinetic equations. II. Convergence proofs for the Boltzmann equation. Comm. Pure Appl. Math. 46 (1993), no. 5, 667-753.

\bibitem{Bardos-Ukai} C. Bardos, S. Ukai,  The classical incompressible Navier-Stokes limit of the Boltzmann equation. Math. Models Methods Appl. Sci. 1 (1991), no. 2, 235-257.

\bibitem{Boltzmann} Boltzmann, L. Weitere Studien \"{u}ber das W\"armegleichgewicht unter Gasmolek\"ulen. Sitzungs. Akad. Wiss. Wien 66 (1872), 275-370; translated as: Further studies on the thermal equilibrium of gas molecules. Kinetic theory, vol. 2, 88-174. Pergamon, London, 1966.

\bibitem{Caflish} R.E.Caflisch, The fluid dynamic limit of the nonlinear Boltzmann equation. Comm.Pure Appl. Math. 33(1980), No.5, 651-666.

\bibitem{Diperna-Lions} R.J. DiPerna, P.-L. Lions, On the Cauchy problem for Boltzmann equations: global existence and weak stability. Ann. of Math. (2) 130 (1989), no. 2, 321-366.

\bibitem{DHWY} R.J. Duan, F.M. Huang, Y. Wang, T. Yang, Global well-posedness of the Boltzmann equation with large amplitude initial data. Arch. Ration. Mech. Anal. 225 (2017), no. 1, 375-424.

\bibitem{E-Guo-M} R. Esposito, Y. Guo,  R. Marra,  Hydrodynamic limit of a kinetic gas flow past an obstacle. Comm. Math. Phys. 364 (2018), no. 2, 765-823.

\bibitem{E-Guo-K-M} R. Esposito, C. Kim, Y. Guo,  R. Marra, Stationary solutions to the Boltzmann equation in the hydrodynamic limit. Ann. PDE 4 (2018), no. 1, pp. 1-119

\bibitem{Golse-Saint-Raymond} F. Golse, L.  Saint-Raymond, The Navier-Stokes limit of the Boltzmann equation for bounded collision kernels. Invent. Math. 155 (2004), no. 1, 81-161.

\bibitem{Guo2006}Y. Guo, Boltzmann diffusive limit beyond the Navier-Stokes approximation. Comm.PureAppl. Math. 59 (2006), no. 5, 626-687.


\bibitem{Guo2010} Y. Guo,  Decay and continuity of the Boltzmann equation in bounded domains. Arch. Ration. Mech. Anal. 197 (2010), no. 3, 713-809.
		
\bibitem{Guo Jang} Y. Guo, J. Jang, Global Hilbert expansion for the Vlasov-Poisson-Boltzmann system. Comm. Math. Phys. 299 (2010), no. 2, 469-501.

\bibitem{Guo Jang Jiang-1} Y. Guo, J. Jang, N. Jiang, Local Hilbert expansion for the Boltzmann equation. Kinet. Relat. Models. 2 (2009), no. 1, 205-214.

\bibitem{Guo Jang Jiang} Y. Guo, J. Jang, N. Jiang, Acoustic limit for the Boltzmann equation in optimal scaling. Comm. Pure Appl. Math. 63 (2010), no. 3, 337-361.



\bibitem{Hilbert} D. Hilbert,  Mathematical problems. Bull. Amer. Math. Soc. (N.S.) 37 (2000), no. 4, 407-436.

\bibitem{Huang-Jiang-Wang} F.M. Huang, S. Jiang, Y. Wang, Zero dissipation limit of full compressible Navier-Stokes equations with a Riemann initial data. Commun. Inf. Syst. 13 (2013), no. 2, 211-246.

\bibitem{Huang-Wang-Yang} F. M. Huang, Y. Wang,  T. Yang, Hydrodynamic limit of the Boltzmann equation with contact discontinuities, Comm. Math. Phys., 295 (2010), pp. 293-326.

\bibitem{Huang-Wang-Yang-1} F. M. Huang, Y. Wang,  T. Yang, Fluid Dynamic Limit to the Riemann Solutions of Euler Equations: I. Superposition of rarefaction waves and contact discontinuity, Kinet. Relat. Models, 3 (2010), pp. 685-728.

\bibitem{Huang-Wang-Wang-Yang} F. M. Huang, Y. Wang, Y. Wang, T. Yang, The limit of the Boltzmann equation to the Euler equations for Riemann problems. SIAM J. Math. Anal. 45 (2013), no. 3, 1741-1811.


\bibitem{Jiang-Masmoudi} N. Jiang, N. Masmoudi, Boundary layers and incompressible Navier-Stokes-Fourier limit of the Boltzmann equation in bounded domain I. Comm. Pure Appl. Math. 70 (2017), no. 1, 90-171.

\bibitem{Lachowicz} M. Lachowicz, On the initial layer and the existence theorem for the nonlinear Boltzmann equation. Math. Methods Appl. Sci. 9 (1987), no. 3, 342-366.

\bibitem{Masmoudi-Raymond}  N. Masmoudi, L. Saint-Raymond,
From the Boltzmann equation to the Stokes-Fourier system in a bounded domain.
Comm. Pure Appl. Math. 56 (2003), no. 9, 1263-1293.

\bibitem{Maxwell} J.C. Maxwell,  On the dynamical theory of gases. Philos. Trans. Roy. Soc. London Ser. A 157 (1867), 49-88. Reprinted in The scientific letters and papers of James Clerk Maxwell, vol. II, 1862-1873, 26-78. Dover, New York, 1965.

\bibitem{Nishida} T. Nishida, Fluid dynamical limit of the nonlinear Boltzmann equation to the level of the compressible Euler equation. Comm. Math. Phys. 61 (1978), no. 2, 119-148.
		
		
\bibitem{Ukai-Asano} S. Ukai, K. Asano, The Euler limit and the initial layer of the nonlinear Boltzmann equation, Hokkaido Math. J., 12 (1983), pp. 303-324.

\bibitem{Xin-1993} Z. P. Xin, Zero dissipation limit to rarefaction waves for the one-dimentional Navier-Stokes equations of compressible isentropic gases, Comm. Pure Appl. Math, 46 (1993), pp. 621-665.

\bibitem{Xin-Zeng} Z. P. Xin and H. H. Zeng, Convergence to the rarefaction waves for the nonlinear Boltzmann equation and compressible Navier-Stokes equations, J. Differential Equations, 249 (2010), pp. 827-871.

\bibitem{Yu}S. H. Yu, Hydrodynamic limits with shock waves of the Boltzmann equations, Comm. Pure Appl. Math, 58 (2005), pp. 409-443.

\end{thebibliography}
\end{document}